# NECESSARY AND SUFFICIENT CONDITIONS OF EXTREMUM FOR POLYNOMIALS AND POWER SERIES IN THE CASE OF TWO VARIABLE

## V. N. Nefedov

*Moscow Aviation Institute (National Research University), Moscow, 125080 Russia*
*e-mail: [nefedovvn54@yandex.ru](mailto:nefedovvn54@yandex.ru)*

The present paper is a continuation of the author's previous works, in which necessary and sufficient local extrema at a stationary point of a polynomial or a power series (and thus of an analytic function) are given. It is known that for the case of one variable, the necessary and sufficient conditions of the extremum are closing, i.e., they can be formulated as a single condition. The next most complicated case is the case with two variables, which is the one considered in this paper. In this case, many procedures, to which the verification of necessary and sufficient conditions is reduced, are based on the computation of real roots of a polynomial from one variable, as well as on the solution of some other rather simple practically realizable problems. An algorithm based on these procedures is described. Nevertheless, there are still cases where this algorithm "doesn't work". For such cases we propose the method of "substitution of polynomials with uncertain coefficients", using which, in particular, we have described an algorithm that allows us to unambiguously answer the question about the presence of a local minimum at a stationary point for a polynomial that is the sum of two $A$-quasi-homogeneous forms, where $A$ - is a two-dimensional vector, whose components are natural numbers.

## 1. INTRODUCTION

The present paper is a continuation of [1], [2], in which some necessary as well as sufficient conditions of extremum (for certainty - minimum) at a stationary point of a polynomial or a degree series (and thus of an analytic function) are given. I will briefly outline the essence of the problem considered in the work. Let $x = (x_1,...,x_m) \in \mathbb{R}^m$, $p(x)$ be a polynomial (or a power series), $0_{(m)} = (0,...,0) \in \mathbb{R}^m$, $p(0_{(m)}) = 0$, $p'(0_{(m)}) = 0_{(m)}$, i.e. $0_{(m)}$ is a stationary point. The classic question from optimization theory is: $0_{(m)}$ is the point of the local minimum or not? Similar problems arise in all branches of mathematics. In such a situation, the first thing that comes to mind is to study the matrix of second derivatives $p''(0_{(m)})$ using the Sylvester criterion. However, this criterion works in such a way that a positive answer to the question of whether a point $0_{(m)}$ is a local minimum is given only if there is a strongly convex function $p(x)$ in the neighborhood $0_{(m)}$ (i.e. its graphical representation is a "bowl"). Otherwise, "more subtle research" is required. And if the matrix of the second derivatives is zero? In a previously published paper [2], it is proposed to use the Newton polyhedron for such studies, i.e. the convex hull of vectors of degrees of variables of polynomial terms. For example, the Newton polyhedron of a polynomial (here $x, y$ are scalar variables)

$$p(x, y) = x^4 y^2 + 2x^2 y^3 + y^4 + 3x^6 y^2 + 3x^4 y^3 + 0.01 x^8 y^3$$

is shown in Fig. 1 (highlighted in dark gray).

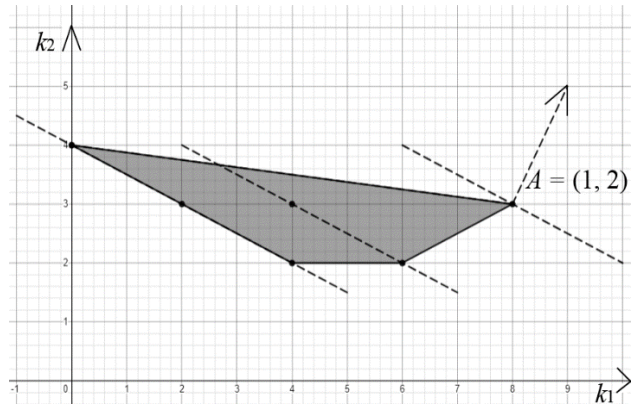

Fig. 1

In [2], some necessary and sufficient conditions are presented and proved based on the study of sums of terms of a polynomial (more generally, a power series) corresponding to the "southwestern" faces of the Newton polyhedron (the polynomial from the above example has three "southwestern" faces: two dimensions 0 (corner points (0, 4), (4, 2)) and one of dimension 1; see Fig. 1). The necessary condition for a local minimum is their non–negativity. Sufficient is non–degeneracy in a weak sense. These results are obtained in [2] for an arbitrary number of variables.

For many cases, this approach already provides a solution to the problem. But there are cases when the necessary condition is fulfilled, but the sufficient one is not. For such cases, a technique using the decomposition of a polynomial (power series) into a sum of quasi-homogeneous forms (analogous to a homogeneous form) is proposed. Figure 1 illustrates the decomposition of a polynomial $p(x, y)$ into the sum of three $A$-quasi-homogeneous forms, where $A$=(1, 2). [2] provides a number of statements regarding such a decomposition with increasingly subtle necessary and sufficient conditions. However, it was not possible to close these conditions for an arbitrary number of variables. Meanwhile, these conditions are closing in case of one variable.

The next in complexity compared to the case of one-variable is to consider the case with two variables, which is what this work is devoted to. For this case, many procedures that are not trivial for the general, turn out to be very simple, easily verifiable in practice. For example, the non-negativity or non-degeneracy of a quasi-homogeneous form is checked by finding the real roots of the polynomial characterizing this form (a simple online procedure). This made it possible to describe simple practically realizable algorithms for solving the problem using the "southwestern" faces of the Newton polyhedron of the polynomial (power series), as well as its decomposition into the sum of $A$-quasi-homogeneous forms, where $A \in \mathbb{N}^2$, where $\mathbb{N}=\{1,2,...\}$ is a natural series.

If these algorithms do not provide an answer to the question being solved, it is proposed to use the method of substitution of polynomials with undefined coefficients. Examples illustrating the operation of the method are given. A statement is proved that allows us to narrow down the infinite set of coefficients for the main terms of polynomials with indeterminate coefficients to two possible variants.

In addition, in this paper, a class of polynomials is distinguished (the sum of two $A$-quasi-homogeneous forms, where $A \in \mathbb{N}^2$) for which a single necessary and sufficient local minimum condition is obtained. On the one hand, this result has an independent meaning, and on the other hand, it allows us to apply sufficiently subtle sufficient local minimum conditions when using statement 15 (the third amplification of theorem 1) from [2], in the formulation of which there is a local minimum condition for the sum of the first few (in particular, two) quasi-homogeneous forms.

It should be noted that the main results of [1], [2], the consequence of which are the main statements of the present paper, are obtained using lemma 7 of [1, p. 208] (see below statement 8, which is its consequence for the two-dimensional case under consideration). The assertion of this lemma is proved in [1], using the method of successive variable elimination from polynomial optimization problems, as well as theorems 21.7, 21.8 of [4], in which a method of finding all solutions of equation $f(x,y)=0$, where $f(x,y)$ is a power series absolutely convergent in the neighborhood of the point $0_{(2)}$, $f(0,0)=0$, and solutions are sought in the form of power series with rational exponents of powers absolutely convergent in the neighborhood of $x=0$, where $y(0)=0$. In [5], the statement of this lemma was generalized from the polynomial case to the case of a power series absolutely convergent in the neighborhood of the zero point.

Wherever it is necessary, the essential difference between the two-dimensional case under consideration and the cases with three and more variables is shown (see remarks 3, 5 below).

## 2. PROBLEM STATEMENT. THE CASE WHEN THE POLYNOMIAL IS A QUASI-HOMOGENEOUS FORM

In this paper we consider polynomials $p(x,y)$ from two real variables $x, y$. Moreover, it will be shown later in Remarks 4, 8 that $p(x,y)$ can be a power series as well. However, for simplicity, for now we speak only about polynomials. Let $p(x,y) \not\equiv 0$, $p(0,0)=0$, $p'(0,0)=(0,\ 0)=0_{(2)}$, i.e. $0_{(2)}$ be a stationary point. The aim of the paper is to obtain practically verifiable necessary as well as sufficient conditions for what $0_{(2)}$ is a point of local minimum $p(x,y)$. Following [2], we will use some notions and notations. We will use the *polynomial carrier* $N_p$ - the set of integer vectors $(\alpha, \beta) \in (\mathbb{N} \cup \{0\})^2$, such that $p(x,y)$ includes a term of the form $ax^\alpha y^\beta$ at $a \neq 0$. The convex hull of this set $\text{Co } N_p$ is called the *Newton polyhedron* of the polynomial $p(x,y)$. If for some vector $A = (A_1, A_2) \in \mathbb{Z}^2 \setminus \{0_{(2)}\}$ there is a number $B \in \mathbb{Z}$, for which the condition $\forall (\alpha, \beta) \in N_p\ \langle A, (\alpha, \beta) \rangle = B$ (where $\langle (x_1, y_1), (x_2, y_2) \rangle = x_1 x_2 + y_1 y_2$ is the scalar product of vectors $(x_1, y_1), (x_1, y_1) \in \mathbb{R}^2$) is satisfied, then the polynomial $p(x,y)$ is called the *$A$-quasi-homogeneous (polynomial) form*. We will call a polynomial $p(x,y)$ a quasi-homogeneous form if for some $A \in \mathbb{Z}^2 \setminus \{0_{(2)}\}$ it is an $A$-quasi-homogeneous form.

**Example 1.** The following polynomials are quasi-homogeneous forms (suitable $A$, $B$ in parentheses):

(а) $p_1(x,y) = 2xy^2 - 3x^2y^3 + 5x^3y^4$ ($A = (1,\ -1)$, $B = -1$);

(б) $p_2(x,y) = 2xy^2 + 3x^3y^2 - 5x^5y^2$ ($A = (0,\ 1)$, $B = 2$);

(в) $p_3(x,y) = 2xy^2 + 3xy^3 + 4xy^5$ ($A = (1,\ 0)$, $B = 1$);

(г) $p_4(x,y) = 2x^4y^2 + 3x^2y^3 + 2y^4$ ($A = (1,\ 2)$, $B = 8$).

Clearly, if a polynomial $p(x,y)$ is an $A$-quasi-homogeneous form, then it is a $\lambda A$-quasi-homogeneous form for any $\lambda \in \mathbb{Q} \setminus \{0\}$, at which $\lambda A \in \mathbb{Z}^2$ (including $\lambda = -1$, where $\mathbb{Q}$ is the set of rational numbers).

Note that if the polynomial $p(x,y)$ is a quasi-homogeneous form, then $\dim \text{Co } N_p < 2$, since for some $A \in \mathbb{Z}^2 \setminus \{0_{(2)}\}$ all points $(\alpha, \beta) \in N_p$ (and hence points from $\text{Co } N_p$) belong to the line given by the equation $A_1(\alpha - \alpha_1) + A_2(\beta - \beta_1) = 0$, where $(\alpha_1, \beta_1)$ is an arbitrary pair from $N_p$. The possible cases are: $\dim \text{Co } N_p = 0$ (the trivial case where $p(x, y)$ consists of a single term), or $\dim \text{Co } N_p = 1$. It is not difficult to show [2] that the converse is also true: if $\dim \text{Co } N_p < 2$, then $p(x, y)$ is a quasi-homogeneous form.

In this section we give a necessary and sufficient condition for what $0_{(2)}$ is the point of local minimum of a polynomial $p(x, y)$ for the simplest case when it is a quasi-homogeneous form. Let us start with the most complicated case when $p(x, y)$ is an $A$-quasi-homogeneous form for some vector $A = (A_1, A_2) \in \mathbb{N}^2$ (further it will be shown that consideration of other cases in which $A \in \mathbb{Z}^2 \setminus (\mathbb{N}^2 \cup \{0_{(2)}\})$ is quite obvious), i.e. has the form $p(x, y) = \varphi_1^A(x, y)$, where for some $s \in \mathbb{N}$ fulfills

$$\varphi_1^A(x, y) = \sum_{i=1}^{s} a_i x^{\alpha_i} y^{\beta_i}, \tag{2.1}$$

$$\begin{aligned} &a_i \neq 0,\ \alpha_i, \beta_i \in \mathbb{N} \cup \{0\},\ A_1\alpha_i + A_2\beta_i = B_1^A \in \mathbb{N}, \\ &i = 1,\ldots,s,\ \alpha_1 > \alpha_2 > \ldots > \alpha_s \geq 0 \end{aligned} \tag{2.2}$$

(herewith $0 \leq \beta_1 < \beta_2 < \ldots < \beta_s$).

For simplicity, we consider that $\text{GCD}(A_1, A_2) = 1$ (otherwise divide $A_1, A_2$ by $\text{GCD}(A_1, A_2)$; $\text{GCD}(A_1, A_2)$ is the greatest common divisor of integers $A_1, A_2$). We will assume that $s \geq 2$, since the case $s = 1$ is obvious. The points $(\alpha_i, \beta_i)$ are on the same line whose coordinates $(\alpha, \beta)$ satisfy the equation $A_1\alpha + A_2\beta = B_1$ (or $\dfrac{\alpha - \alpha_1}{-A_2} = \dfrac{\beta - \beta_1}{A_1}$), with a guide vector $e = (-A_2, A_1)$, and hence $(\alpha_i, \beta_i) = (\alpha_1, \beta_1) + v_i e$, $i = 1, \ldots, s$, where $v_1 = 0$, $v_i > 0$, $i = 2, \ldots, s$. Using the fact that $\text{GCD}(|e_1|, |e_2|) = \text{GCD}(A_1, A_2) = 1$, $(\alpha_i - \alpha_1, \beta_i - \beta_1) = v_i e$, we obtain that $v_i \in \mathbb{N}$, $i = 2, \ldots, s$, $v_1 = 0$. Thus, at $x \neq 0$ we have:

$$\begin{aligned} \varphi_1^A(x, y) &= \sum_{i=1}^{s} a_i x^{\alpha_i} y^{\beta_i} = \sum_{i=1}^{s} a_i x^{\alpha_1 + v_i e_1} y^{\beta_1 + v_i e_2} = x^{\alpha_1} y^{\beta_1} \sum_{i=1}^{s} a_i \left(x^{e_1} y^{e_2}\right)^{v_i} = \\ &= x^{\alpha_1} y^{\beta_1} \sum_{i=1}^{s} a_i u^{v_i} = x^{\alpha_1} y^{\beta_1} g_1^A(u),\ u = x^{e_1} y^{e_2},\ e = (-A_2, A_1),\ g_1^A(u) = \sum_{i=1}^{s} a_i u^{v_i}, \end{aligned} \tag{2.3}$$

where $g_1^A(u)$ is a polynomial from one variable $u \in \mathbb{R}$. In this case, by virtue of $\text{GCD}(|e_1|, |e_2|) = 1$ at least one of the numbers among $e_1, e_2$, is odd, and hence the variable $u = x^{e_1} y^{e_2}$ can take on $(x, y) \in \mathbb{R}^2$ any real values. Hereafter we will call $g_1^A(u)$ the *characteristic polynomial* for the quasi-homogeneous form $\varphi_1^A(x, y)$, and $a_1 x^{\alpha_1} y^{\beta_1}$ its *main term*. Accordingly, in the case $x = 0$ we have: $p(0, y) = \varphi_1^A(0, y) \equiv 0$ at $\alpha_s > 0$, and $p(0, y) = a_s y^{\beta_s}$ at $\alpha_s = 0$. Note that the term "characteristic polynomial" is used in matrix theory, while here this term is used for quasi-homogeneous forms and has nothing to do with matrices.

The following is true

**Assertion 1.** A polynomial $p(x, y) = \varphi_1^A(x, y)$ of the form (2.1) - (2.3), where $A = (A_1, A_2) \in \mathbb{N}^2$, is non-negative if and only if the following conditions are satisfied:

1) $a_1 > 0$, $a_s > 0$, $\alpha_1, \beta_1, \alpha_s, \beta_s \in 2\mathbb{N} \cup \{0\}$;

2) polynomial $g_1^A(u) = \sum_{i=1}^{s} a_i u^{v_i}$ is non-negative, i.e. (by virtue of the fact that $g_1^A(0) = a_1 > 0$; see the previous condition) either has no real roots, or all its real roots have even multiplicity.

**Proof.** Let us prove necessity (sufficiency is obvious). Let us first show that 1) is true. For example, if $a_1 > 0$, $\alpha_1$ is odd, $\beta_1$ is any, then $p(-t^{A_1}, t^{A_2+1}) = -a_1 t^{B_1^A + \beta_1} + o(t^{B_1^A + \beta_1})$, i.e., the polynomial $p(x, y)$ is not non-negative. Accordingly, if $a_1 > 0$, $\alpha_1$ is even, $\beta_1$ is odd, i.e. $p(t^{A_1}, -t^{A_2+1}) = -a_1 t^{B_1^A + \beta_1} + o(t^{B_1^A + \beta_1})$, the polynomial $p(x, y)$ is not non-negative again. The other cases of non-fulfillment of condition 1) are treated similarly. Let us now show the validity of 2). Suppose that $g_1^A(u_0) < 0$ is satisfied for some $u_0 \in \mathbb{R}$. Then $u_0 \neq 0$ (since $g_1^A(0) = a_1 > 0$). Using the fact that at least one of the numbers among $e_1, e_2$, is odd, it is easy to find $x \neq 0, y \neq 0$, which are solutions of the equation $u_0 = x^{e_1} y^{e_2}$ (for example, if $e_1$ is odd, we assume: $y = 1, x = (u_0)^{1/e_1}$). Then by virtue of (2.3) for the chosen $x, y$ ones it is fulfilled (taking into account parity $\alpha_1, \beta_1$) $p(x, y) = x^{\alpha_1} y^{\beta_1} g_1^A(u_0) < 0$.

Thus, we obtain a necessary and sufficient condition for non-negativity of a polynomial $p(x, y) = \varphi_1^A(x, y)$ of the form (2.1) - (2.3), where $A = (A_1, A_2) \in \mathbb{N}^2$. Obviously, in the case of non-negativity of this polynomial $0_{(2)}$ is the point of its local (and even global) minimum, and in the case where non-negativity is not satisfied, it is not the point of the local minimum of this polynomial. Indeed, if we find a point $(x^*, y^*) \in \mathbb{R}^2$, for which $\varphi_1^A(x^*, y^*) < 0$, then for $u_1(t) = x^* t^{A_1}$, $u_2(t) = y^* t^{A_2}$ for all $t > 0$ (including arbitrarily small) we have:
$\varphi_1^A(u_1(t), u_2(t)) = \varphi_1^A(x^*, y^*) t^{B_1^A} < 0$.

**Example 2.** Note that the polynomial $p_4(x, y)$ from Example 1 is an $A$-quasi-homogeneous form at $A = (1, 2)$, i.e., it is in the scope of Assertion 1. For this polynomial, condition 1) is satisfied, the characteristic polynomial $g_1^A(u)$, satisfying (2.3), is of the form $g_1^A(u) = 2 + 3u + 2u^2$, where $u = x^{-2} y$. This square trinomial has a negative discriminant, i.e., it has no real roots, and hence, by virtue of statement 1, $0_{(2)}$ is the point of the local (and global) minimum of the polynomial $p_4(x, y)$.

Consider now the case when a polynomial $p(x, y) = \varphi_1^A(x, y)$ of the form (2.1) - (2.3), is an $A$-quasi-homogeneous form for some $A = (A_1, A_2) \in \mathbb{Z}^2 \setminus \{0_{(2)}\}$, such that $A \notin \mathbb{N}^2$. Since, as already noted, any $A$-quasi-homogeneous form where $A \in \mathbb{Z}^2 \setminus \{0_{(2)}\}$, will be simultaneously a $(-A)$-quasi-homogeneous form, it suffices to restrict ourselves to considering the cases (a) $A_1 > 0$, $A_2 < 0$, (b) $A_1 = 0$, $A_2 > 0$, (c) $A_1 > 0$, $A_2 = 0$ (e.g., if $A_1 < 0, A_2 > 0$, then $-A_1 > 0$, $-A_2 < 0$), and at the same time with $A \notin \mathbb{N}^2$ we assume that $-A \notin \mathbb{N}^2$.

Let us make some changes in condition (2.2). In cases (a), (b) we replace the condition $\alpha_1 > \alpha_2 > ... > \alpha_s \geq 0$ in (2.2) by $0 \leq \alpha_1 < \alpha_2 < ... < \alpha_s$ (then in case (a) $0 \leq \beta_1 < \beta_2 < ... < \beta_s$, and in case (b) $\beta_1 = \beta_2 = ... = \beta_s$; see polynomials $p_1(x, y)$, $p_2(x, y)$, in Example 1). Accordingly, in case (c) we replace the condition $\alpha_1 > \alpha_2 > ... > \alpha_s \geq 0$ with $0 \leq \beta_1 < \beta_2 < ... < \beta_s$ (in this case, $\alpha_1 = \alpha_2 = ... = \alpha_s$; see the polynomial $p_3(x, y)$ in Example 1). In each of these cases $a_1 x^{\alpha_1} y^{\beta_1}$ is the main term of the polynomial $p(x, y) = \varphi_1^A(x, y)$ (any other term of this polynomial has a vector of degrees greater than the Pareto vector of $(\alpha_1, \beta_1)$), which determines whether $0_{(2)}$ will be a point of local minimum of this polynomial or not. For $0_{(2)}$ to be a point of local minimum, it is necessary and sufficient for this term to be non-negative, i.e., at $a_1 > 0$, $\alpha_1, \beta_1 \in 2\mathbb{N} \cup \{0\}$. In Example 1, for all three polynomials, this term is $2xy^2$, which can take negative values in any small neighborhood of the point $0_{(2)}$, i.e., it is not the point of local minimum of these polynomials. Indeed, for each of these polynomials, substituting $x = -t$, $y = t$ gives a polynomial that is negative in a sufficiently small neighborhood of the point $t = 0$ at $t > 0$. Quite similarly and in the general case, if the condition $a_1 > 0$, $\alpha_1, \beta_1 \in 2\mathbb{N} \cup \{0\}$ is not satisfied, it is not difficult to choose numbers $C_1, C_2 \in \{1, -1\}$ for which $a_1 C_1^{\alpha_1} C_2^{\beta_1} < 0$, and thus in the neighborhood of the point $t = 0$ for sufficiently small $t > 0$ the condition $\varphi_1^A(C_1 t, C_2 t) = a_1 C_1^{\alpha_1} C_2^{\beta_1} t^{\alpha_1 + \beta_1} + o(t^{\alpha_1 + \beta_1}) < 0$ is true.

**Remark 1.** It was noted earlier that the choice of the vector $A \in \mathbb{Z}^2 \setminus \{0_{(2)}\}$ for the $A$-quasi-homogeneous form is not unambiguous. However, by imposing some simple additional conditions on $A$, one can arrive at unambiguity of this choice. For example, for a quasi-homogeneous form consisting of a single term, it can be chosen arbitrarily (in particular, $A = (1, 1)$). In this case, the verification of the non-negativity of this polynomial is trivial (see condition 1) of statement 1). Let now the quasi-homogeneous form contains at least two terms $(\alpha_1, \beta_1)$, $(\alpha_2, \beta_2)$. Then, if $A_1 \neq 0$, $A_2 \neq 0$, then the additional condition $A_1 > 0$ leads to the system
$$A_1 > 0, \quad A_1(\alpha_2 - \alpha_1) + A_2(\beta_2 - \beta_1) = 0,$$
which, when adding the condition $\text{GCD}(A_1, |A_2|) = 1$, always has a single solution, which we denote by $\bar{A}$. It is clear that if we choose another vector satisfying this system without additional condition $\text{GCD}(A_1, |A_2|) = 1$, we will get a vector $A = \lambda \bar{A}$ at $\lambda \in \mathbb{N}$ (we take into account the conditions: $A_1 > 0$, $\text{GCD}(A_1, |A_2|) = 1$). With this change, the result of investigating the point $0_{(2)}$ for a local minimum will remain the same, i.e., we can always restrict ourselves to using $A = \bar{A}$. Accordingly, in the case of $A_1 = 0$ it is sufficient to take $A = (0, 1)$, and in the case of $A_1 = 0$ - $A = (1, 0)$.

**Remark 2**. Simple examples show that the above statements for the case of two variables do not carry over to the case of three variables, indicating that this case (as well as the case with one variable) is exceptional. For example, the polynomial $p(x, y, z) = x^2 + 2xy + y^2 + x^2 y^2 z$ is a $A$-quasi-homogeneous form where $A = (1, 1, -2) \in \mathbb{Z}^3 \setminus \{0_{(3)}\}$, $A \notin \mathbb{N}^3$. This polynomial has not one but three main terms, which are not Pareto-comparable. The sum of these terms gives a non-negative $(1, 1)$-quasi-homogeneous form of two variables. In this case $0_{(3)}$ is not a point of local minimum of the original polynomial, since $p(t, -t, -t) = -t^5 < 0$ for all $t > 0$.

3. CASE WHEN NEWTON'S POLYHEDRON HAS DIMENSION 2. USE OF MAIN QUASI-HOMOGENEOUS FORMS

Let now $\dim \text{Co} N_p = 2$. In this case the polynomial $p(x, y)$ cannot be a quasi-homogeneous form (see Section 2). Nevertheless, it can have the main quasi-homogeneous forms corresponding to the faces of the Newton polyhedron $\text{Co} N_p$. Each eigenface $C \subset \text{Co} N_p$ corresponds to a main quasi-homogeneous form $\varphi_C(x, y) = p_{N_p \cap C}(x, y)$, where for any nonempty set $N \subseteq N_p$, $p_N(x, y)$ is the sum of the terms of the polynomial $p(x, y)$ whose degree vectors belong to $N$ (the $N$-*shortening* of the polynomial $p(x,y)$). Thus, $N_{\varphi_C} = N_p \cap C$ and all members of the polynomial $\varphi_C(x,y)$ are members of the polynomial $p(x,y)$.

It should be noted that Newton's polyhedron is a tool for investigating a wide class of problems (see, e.g., [6]-[9]). In particular, Newton's theory of polyhedra connects the geometry of polyhedra with algebraic geometry [9].

Here we may recall that a convex non-empty set $C \subseteq \mathbb{R}^n$ is called a *face* of a convex closed set $Y \subseteq \mathbb{R}^n$, if $C \subseteq Y$ and $\forall y^{(1)}, y^{(2)} \in Y$, $\forall \alpha \in (0, 1)$ in the case of $\alpha y^{(1)} + (1-\alpha) y^{(2)} \in C$ the condition $y^{(1)}, y^{(2)} \in C$ is true. A face $C$ is called an *eigenface* if $C \neq Y$. A point $y \in Y$ is called *corner (angular)* if $\{y\}$ is a face.

In the considered case with $\dim \operatorname{Co} N_p = 2$, the main quasi-homogeneous forms can correspond only to the proper faces of the polyhedron $\operatorname{Co} N_p$: corner points (faces of dimension 0) and sides of $\operatorname{Co} N_p$ (faces of dimension 1).

It should be noted that there is also another definition of the main quasi-homogeneous form of the polynomial [2]. A polynomial $\varphi^A(x,y)$ is called the main $A$-quasi-homogeneous form of a polynomial $p(x,y)$, where $A = (A_1, A_2) \in \mathbb{Z}^2 \setminus \{0_{(2)}\}$, if all terms of $\varphi^A(x,y)$ are members of $p(x,y)$ and $N_{\varphi^A} = \operatorname{Arg\,min}\{\langle A, k \rangle | k \in N_p\}$. A polynomial $\varphi(x,y)$ is called the main quasi-homogeneous form of a polynomial $p(x,y)$ if for some $A \in \mathbb{Z}^2 \setminus \{0_{(2)}\}$ it is its main $A$-quasi-homogeneous form. In [2] it is shown that both definitions are equivalent. In particular, this means that for any eigenface $C$ of a polyhedron $\operatorname{Co} N_p$ there is a vector $A \in \mathbb{Z}^2 \setminus \{0_{(2)}\}$, such that $p_{N_p \cap C}(x,y)$ is the main $A$-quasi-homogeneous form of the polynomial $p(x,y)$, i.e., $N_p \cap C = \operatorname{Arg\,min}\{\langle A, k \rangle | k \in N_p\}$. As will be seen later, finding a suitable vector $A$ for the faces of the polyhedron $\operatorname{Co} N_p$ considered below is not a difficult task.

Accordingly, if for some $A \in \mathbb{Z}^2 \setminus \{0_{(2)}\}$ $\varphi^A(x,y)$ is the main $A$-quasi-homogeneous form of a polynomial $p(x,y)$, whence $N_{\varphi^A} = \operatorname{Arg\,min}\{\langle A, k \rangle | k \in N_p\}$, then, as shown in [2], it corresponds to an eigenface (of a polyhedron $\operatorname{Co} N_p$) $C_A = \operatorname{Arg\,min}\{\langle A, k \rangle | k \in \operatorname{Co} N_p\}$, and it follows that: $N_{\varphi^A} = N_p \cap C_A$, $C_A = \operatorname{Co}(N_p \cap C_A) = \operatorname{Co} N_{\varphi^A}$, $\varphi^A(x,y) = p_{N_p \cap C_A}(x,y)$.

Using only the main quasi-homogeneous forms of the polynomial, we can already obtain some necessary and sufficient conditions for the local minimum of the polynomial $p(x,y)$ at the point $0_{(2)}$. We will need the following definitions. Let a function $f(x,y)$ be defined on the entire space $\mathbb{R}^2$. Say that a function $f(x,y)$ is *non-negative* if $\forall (x,y) \in \mathbb{R}^2$ $f(x,y) \geq 0$, a function $f(x,y)$ is *nondegenerate in the weak sense* if $[x \neq 0, y \neq 0] \Rightarrow f(x,y) \neq 0$ is satisfied. In this section we will use the following corollaries of Theorem 2, Lemma 8, and Remark 3 of [1].

**Theorem 1** [1, p. 209]. Let $p(x,y)$ be a polynomial, $p(0,0) = 0$, $p'(0,0) = 0_{(2)}$, and for all $A \in \mathbb{N}^2$ the main $A$-quasi-homogeneous forms of the polynomial $p(x,y)$ are non-negative and nondegenerate in the weak sense. Then $0_{(2)}$ is the point of local minimum of $p(x,y)$.

**Lemma 1** [1, p. 210]. Let $0_{(2)}$ be the point of local minimum of the polynomial $p(x,y)$, $p(0,0) = 0$, $p'(0,0) = 0_{(2)}$. Then for all $A \in \mathbb{N}^2$ the main $A$-quasi-homogeneous forms of the polynomial $p(x,y)$ are non-negative.

**Remark 3.** It should be noted that the statements of Theorem 1 and Lemma 1 do not consider the main $A$-quasi-homogeneous forms of the polynomial $p(x, y)$ for any vector $A \in \mathbb{Z}^2 \setminus \{0_{(2)}\}$, but only for $A \in \mathbb{N}^2$. Meanwhile, in many other problems may be required the main $A$-quasi-homogeneous of the polynomial $p(x, y)$ for all vectors $A \in \mathbb{Z}^2 \setminus \{0_{(2)}\}$, e.g., when investigating the polynomial for strict positivity in positive orthant (see [11], as well as references to other works cited there). For example, as shown in [2], when investigating a quasi-homogenious polynomial form depending on more than two variables for non-negativity and nondegeneracy in the weak sense, it may be necessary to check the set of polynomials generated by this form for strict positivity in the positive orthant (see, moreover, Remark 5 below).

Note that in the case $A \in \mathbb{N}^2$ for the main $A$-quasi-homogeneous form $\varphi^A(x, y) = p_{N_{\varphi^A}}(x, y)$, where $N_{\varphi^A} = \operatorname{Arg\,min}\{\langle A, k \rangle | k \in N_p\}$, and its corresponding face $C_A = \operatorname{Co} N_{\varphi^A}$ (for which $N_{\varphi^A} = N_p \cap C_A$) holds

$$N_{\varphi^A} = \operatorname{Arg\,min}\{\langle A, k \rangle | k \in N_p\} \subseteq P(N_p),\ C_A \subseteq P(\operatorname{Co} N_p),$$

where $P(Y) = \{y^* \in Y\,|\,\forall y \in Y\ [y_1 \leq y_1^*, y_2 \leq y_2^*] \Rightarrow y = y^*\}$ is the set of Pareto-optimal points of an arbitrary set $Y \subset \mathbb{R}^2$. A simple statement has been used here:

$$C_A = \operatorname{Arg\,min}\{\langle A, k \rangle | k \in \operatorname{Co} N_p\} \subseteq P(\operatorname{Co} N_p).$$

Thus, using Theorem 1 and Lemma 1, we will need main quasi-homogeneous forms corresponding not to all faces, but only to the faces contained in the contained in the contained in the "southwestern" part of the boundary of polyhedron $\operatorname{Co} N_p$ ("southwest" faces). These faces are easily distinguished even visually. In addition, [2], [10] describe practically realizable algorithms for selecting suitable $A \in \mathbb{N}^2$ to isolate all main $A$-quasi-homogeneous forms of the polynomial corresponding to the cases $A \in \mathbb{N}^2$.

**Remark 4.** In the case of $A \in \mathbb{N}^2$, we can set the problem of finding the set of main $A$-quasi-homogeneous forms not only for a polynomial, but also for a power series, which will be the result of the expansion of the analytic function by powers of variables in the neighborhood of a stationary point, since in this case their number is finite and a finite number of "first" terms from this expansion will be needed. In this case, Theorem 1 and Lemma 1 remain valid for the power series $p(x, y)$ [5].

Note now that the set of all main quasi-homogenious forms of a polynomial $p(x, y)$ (or of a power series; see Remark 4) which are main $A$-quasi-homogeneous at some $A \in \mathbb{N}^2$ can be divided into three groups.

*Group 1.* This group includes the main quasi-homogeneous forms of the polynomial $p(x, y)$, corresponding to the faces of the polyhedron $\operatorname{Co} N_p$ of dimension 0. Each of these forms is a single term of the form $ax^\alpha y^\beta$, which is a term of the polynomial $p(x, y)$. Thus we write $coef(p, (\alpha, \beta)) = a$. As shown in [2], these single terms correspond to the corner points of the polyhedron $\operatorname{Co} N_p$, i.e., $(\alpha, \beta) \in \Psi(N_p) \subseteq N_p$, where through $\Psi(N_p)$ denotes in [2] the set of corner points of $\operatorname{Co} N_p$. Moreover, the vectors $(\alpha, \beta)$ are Pareto optimal on the set $N_p$, i.e., $(\alpha, \beta) \in P(N_p)$. Thus, $(\alpha, \beta) \in \Omega(N_p)$, where according to the notation from [2] $\Omega(N_p) = P(N_p) \cap \Psi(N_p)$. The converse of [2] is also true: the single term of the polynomial $p(x, y)$, corresponding to each point in $\Omega(N_p)$, is for some $A \in \mathbb{N}^2$ a main $A$-quasi-homogenious form of this polynomial. The extraction (even visually, using a graphical representation of the polyhedron $\operatorname{Co} N_p$) of points from $\Omega(N_p)$ is a rather simple task (see also [2], [10] for algorithms of their extraction for an arbitrary number of variables). A necessary and sufficient condition for non-negativity of such forms, which in this case also entails nondegeneracy in the weak sense, is $a > 0, \alpha, \beta \in 2\mathbb{N} \cup \{0\}$ (see condition 1) of the statement 1).
Thus, a necessary condition for what $0_{(2)}$ is a point of local minimum of a polynomial $p(x, y)$ is

$$\forall (\alpha, \beta) \in \Omega(N_p) \quad coef(p, (\alpha, \beta)) > 0, \ \alpha, \beta \in 2\mathbb{N} \cup \{0\}. \tag{3.1}$$

We will further assume that this condition is satisfied (in the first case of its non-fulfillment the problem under consideration is solved and $0_{(2)}$ is not a point of local minimum of the polynomial $p(x, y)$).

*Group 2.* This group will include the main quasi-homogeneous forms of the polynomial $p(x, y)$ of dimension 1 (i.e., corresponding to the sides of the polyhedron $Co N_p$), which are the sum of two single terms, each of which in this case corresponds to one of the two corner points of that side, which are corner points of $Co N_p$ and simultaneously belong to $P(N_p)$ (i.e., each of them belongs to $\Omega(N_p)$ and thus belongs to group 1). Thus, group 2 will include some sums of two quasi-homogenious forms from group 1, which are single terms. Since the condition (3.1) is supposed to be fulfilled for the single terms included in this sum, which ensures non-negativity and nondegeneracy in the weak sense of each of them, the sum of these single terms will obviously also be non-negative and nondegenerate in the weak sense.

*Group 3.* This group will include all other main quasi-homogeneous forms of the polynomial $p(x, y)$ (not included in groups 1, 2) corresponding to the faces (sides) of the polyhedron $Co N_p$, of dimension 1, which are sums of at least three single terms: two single terms corresponding to the corner points of this face (i.e., included in group 1), and there is at least one single term corresponding to an intermediate point of this face located between its corner points.

Thus, in the case when condition (3.1) is satisfied and all main $A$-quasi-homogeneous forms of the polynomial $p(x, y)$, corresponding to cases $A \in \mathbb{N}^2$, belong to groups 1 and 2, then by virtue of Theorem 1 $0_{(2)}$ is the point of local minimum of the polynomial $p(x, y)$. Let now the set of quasi-homogeneous forms from group 3 is not empty and one of such forms has been singled out. Then for some $A \in \mathbb{N}^2$ it will be a polynomial $\varphi_1^A(x, y)$ of the form (2.1), (2.2).

The problem of determining the vector $A = (A_1, A_2) \in \mathbb{N}^2$ from the terms of the extracted quasi-homogeneous form (e.g., they can be determined visually from the image $Co N_p$; see also in [2], [10] for algorithms of their extraction for an arbitrary number of variables) is a simple computational problem. This vector is determined by any two terms included in this form, e.g., $a_1 x^{\alpha_1} y^{\beta_1}$, $a_2 x^{\alpha_2} y^{\beta_2}$ from the equation $A_1 \alpha_1 + A_2 \beta_1 = A_1 \alpha_2 + A_2 \beta_2$, whence $\frac{A_1}{A_2} = \frac{\beta_2 - \beta_1}{\alpha_1 - \alpha_2}$. We reduce the right-hand side of this equality to an irreducible fraction of the form $\frac{m}{n}$, where $m, n \in \mathbb{N}$, and then suppose $A_1 = m, A_2 = n$. Thus $GCD(A_1, A_2) = 1$.

We denote the set of vectors $A = (A_1, A_2) \in \mathbb{N}^2$ satisfying condition $GCD(A_1, A_2) = 1$ by $\mathbb{N}_0^2$.

At this point, using lemma 1, we can use the statement 1 obtained earlier. By virtue of lemma 1, the non-negativity of $\varphi_1^A(x, y)$ is only a necessary condition that $0_{(2)}$ is the point of local minimum of the polynomial $p(x, y)$. Therefore, one can also use assertion 1 only to obtain a negative result. Namely, if the main $A$-quasi-homogeneous form $\varphi_1^A(x, y)$ of the polynomial $p(x, y)$ is not non-negative (i.e., it does not satisfy the necessary and sufficient conditions from statement 1), then $0_{(2)}$ is not a point of local minimum of the polynomial $p(x, y)$.

If all the main $A$-quasi-homogeneous forms of a polynomial $p(x, y)$ (or of a power series; see Remark 4) from group 3 are non-negative, then the necessary condition from Lemma 1 is satisfied and we proceed to check the sufficient conditions. The strongest of them is the check of quasi-homogeneous forms from group 3 for nondegeneracy in the weak sense. If all of them turn out to be nondegenerate in the weak sense, then according to Theorem 1 the point $0_{(2)}$ is the point of the local minimum of the polynomial $p(x, y)$. Here we can use

**Assertion 2.** A quasi-homogeneous form $\varphi_1^A(x, y)$ of the form (2.1) - (2.3) satisfying condition 1) of statement 1 is nondegenerate in the weak sense if and only if $g_1^A(u) > 0$ at all $u \in \mathbb{R}$, i.e. (by virtue of the fact that $g_1^A(0) = a_1 > 0$) this polynomial has no real roots.

**Proof. Sufficiency.** Let $g_1^A(u) > 0$ be satisfied for all $u \in \mathbb{R}$. By Proposition 1, the quasi-homogeneous form $\varphi_1^A(x, y)$ is non-negative. Let us show that it is nondegenerate in the weak sense. Suppose that we find a point $(x_0, y_0) \in \mathbb{R}^2$ such that $x_0 \neq 0$, $y_0 \neq 0$, $\varphi_1^A(x_0, y_0) = 0$. Then for $u_0 = x_0^{e_1} y_0^{e_2}$, where $e = (-A_2, A_1)$, by virtue of (2.3) the equality $g_1^A(u_0) = 0$, which contradicts the initial assumption, is satisfied.

**Necessity.** Suppose that the quasi-homogeneous form $\varphi_1^A(x, y)$ is nondegenerate in the weak sense. Let us show that $g_1^A(u) > 0$ for all $u \in \mathbb{R}$. Suppose that we find a point $u_0 \in \mathbb{R}$ such that $g_1^A(u_0) = 0$. Since $g_1^A(0) = a_1 > 0$, then $u_0 \neq 0$. Using the fact that at least one of the numbers among $e_1, e_2$ is odd (due to the fact that $GCD(A_1, A_2) = 1$), it is easy to choose $x_0 \neq 0$, $y_0 \neq 0$, which are solutions of the equation $u_0 = x_0^{e_1} y_0^{e_2}$ (for example, if $e_1$ is odd, we assume: $y_0 = 1, x_0 = (u_0)^{1/e_1}$). Then, by virtue of (2.3), $\varphi_1^A(x_0, y_0) = x_0^{\alpha_1} y_0^{\beta_1} g_1^A(u_0) = 0$ is satisfied for the chosen $x_0, y_0, u_0$, which contradicts the nondegeneracy of $\varphi_1^A(x, y)$ in the weak sense.

Thus, to check the main quasi-homogeneous forms of group 3 for nondegeneracy in the weak sense, it is enough to find the real roots of their characteristic polynomials. If all of them turn out to be nondegenerate in the weak sense (i.e., there are no real roots), then, according to Theorem 1, the point $0_{(2)}$ is the point of local minimum of the polynomial (or degree series; see Remark 4) of $p(x, y)$. In the case that for some $A = (A_1, A_2) \in \mathbb{N}_0^2$ we find a main quasi-homogeneous form $\varphi_1^A(x, y)$ from group 3 which is non-negative but for which the condition of nondegeneracy in the weak sense is not satisfied, i.e., its characteristic polynomial has real roots (at least one) of even multiplicity, then more subtle investigations will be required, given below. Such a form from group 3 may not be the only one. Other forms will correspond to other $A = (A_1, A_2) \in \mathbb{N}_0^2$ and each of them should be considered separately.

It should be noted that in the considered two-dimensional case the method considered in this section is ideologically close to the Newton diagram method. For example, the main quasi-homogeneous forms from groups 2 and 3 are analogs of the *reference polynomials* [4] in the Newton diagram method.

**Example 3.** Let $p_a(x, y) = x^2 y^6 - (2+a)x^4 y^5 + x^6 y^4 + y^{10} - 10xy^9 - 0.1x^8 y^4$
- a polynomial, where $a$ is a parameter, $a \neq -2$. Then
$$N_{p_a} = \{(2, 6), (4, 5), (6, 4), (0, 10), (1, 9), (8, 4)\}$$

(see the image of $N_{p_a}$, Co $N_{p_a}$ in Fig. 2). For this polynomial, condition (3.1) is satisfied (here $\Omega(N_{p_a}) = \{(0, 10), (2, 6), (6, 4)\}$) and the only quasi-homogeneous form from group 3 is $\varphi_a(x, y) = x^6 y^4 - (2+a)x^4 y^5 + x^2 y^6$ (it is a (1, 2)-quasi-homogeneous form) whose characteristic polynomial is $g_a(u) = 1 - (2+a)u + u^2$. Then at $a = 0.01$ the polynomial $g_a(u) = g_{0.01}(u)$ is not non-negative ($g_{0.01}(1) = -0.01 < 0$), and hence by virtue of statement 1 the quasi-homogeneous form $\varphi_{0.01}(x, y)$ is not non-negative, and then by virtue of lemma 1 $0_{(2)}$ is not the point of the local minimum of the polynomial $p_{0.01}(x, y)$. At $a = -0.01$, the polynomial $g_a(u) = g_{-0.01}(u)$ has no real roots since $g_{-0.01}(u) = 1 - 1.99u + u^2 = 0.995(1-u)^2 + 0.005(1+u^2) > 0$ at all $u \in \mathbb{R}$. But then by statement 2 the quasi-homogeneous form $\varphi_{-0.01}(x, y)$ is non-negative and nondegenerate in the weak sense, and hence by Theorem 1 $0_{(2)}$ is the point of local minimum of the polynomial $p_{-0.01}(x, y)$. At $a = 0$, the polynomial $g_a(u) = g_0(u) = 1 - 2u + u^2 = (1-u)^2$ is non-negative and has a real root $u = 1$ of multiplicity 2. Thus, in this case the conditions of Lemma 1 are not violated, and Theorem 1 "does not work", and hence further investigation will require more subtle methods, which will be given below.

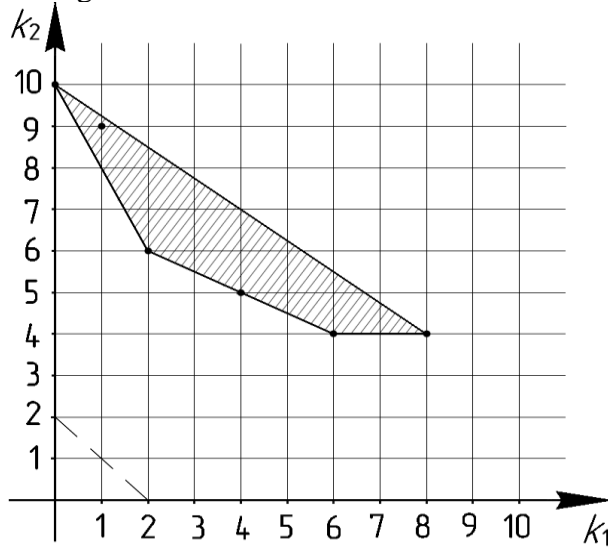

Fig. 2

**Remark 5.** Let us again emphasize the essential difference between the case of two and a large number of variables. The statements of Lemma 1 and Theorem 1 are valid for an arbitrary number of variables [2]. However, the verification of non-negativity and nondegeneracy in the weak sense of quasi-homogenious forms essentially depending on three or more variables will require the study of a finite set of polynomials with rational exponents of powers of the variables for non-negativity and strict positivity in the positive orthant [2]. Such a problem is investigated in [11], where is proposed an algorithm for its solution, which is significantly more complicated than the computation of real roots from a polynomial (to which the case of two variables is reduced). It may be noted here that there is a possibility to simplify the problem by replacing the check for a quasi-homogenious form of nondegeneracy in the weak sense by nondegeneracy in the strong sense [2] (i.e. to a stronger sufficient condition), which in the case of non-negativity of this form is satisfied if and only if the quasi-homogenious form depending on the variables $x_1,...,x_m$, is strictly positive on the surface of the unit cube

$Q = \{(x_1,...,x_m) \in \mathbb{R}^m \mid |x_i| \leq 1, i = 1,...,m\}$. The latter problem reduces to solving several global optimization problems of a polynomial (i.e., a Lipschitz function) on the unit cube (see some methods of solving this problem in [12] and other works of these authors).

# 4. USING THE POLYNOMIAL DECOMPOSITION INTO THE SUM OF $A$-QUASI-HOMOGENEOUS FORMS

In this section we will assume that the Newton polyhedron of a polynomial (or of a degree series; see Remark 8 below) $p(x, y)$ again has dimension 2 and that all main $A$-quasi-quasi-homogenious forms of this polynomial are non-negative at all $A \in \mathbb{N}^2$, i.e. in particular condition (3.1) is satisfied. In this connection, let us introduce into consideration for a given polynomial $p(x, y)$ the set $A_p \subset \mathbb{N}_0^2$ of integer vectors $A \in \mathbb{N}_0^2$ such that $A \in A_p$ if and only if it is simultaneously satisfied that:

1) the main $A$-quasi-homogenious form $\varphi_1^A(x, y)$ of the polynomial $p(x, y)$ belongs to group 3 (i.e., it is the sum of at least three single terms);

2) the characteristic polynomial $g_1^A(u)$, defined according to (2.3), is non-negative (see statement 1);

3) the polynomial $g_1^A(u)$ has real roots.

Let us denote $U_p(A) = \{u \in \mathbb{R} \mid g_1^A(u) = 0\}$. Note that since $g_1^A(0) = a_1 > 0$, then $0 \notin U_p(A)$.

If $A_p = \varnothing$, then, as we have seen, the question of whether $0_{(2)}$ is a point of local minimum of a polynomial $p(x, y)$ is simply solved using Statement 2 and Theorem 1, i.e., in this case $0_{(2)}$ is a point of local minimum of this polynomial.

Let now be $A_p \neq \varnothing$. Then for each $A \in A_p$ we consider the decomposition of the polynomial $p(x, y)$ into a sum of $A$-quasi-homogenious forms (such a decomposition can be obtained for any $A \in \mathbb{Z}^2 \setminus \{0_{(2)}\}$):

$$p(x, y) = \varphi_1^A(x, y) + \varphi_2^A(x, y) + ... + \varphi_{r_A}^A(x, y), \qquad (4.1)$$

$$\varphi_i^A(x, y) \not\equiv 0; \forall k \in N_{\varphi_i^A} \langle A, k \rangle = B_i^A \in \mathbb{N}, \ i = 1, ..., r_A, \ B_1^A < B_2^A < ... < B_{r_A}^A, \qquad (4.2)$$

where $r_A \geq 2$ (since the case of $\dim \mathrm{CoN}_p = 2$) is considered.

Let $H_{\varphi_i^A} = \{(x, y) \in \mathbb{R}^2 \mid \varphi_i^A(x, y) = 0\}$, $i = 1, ..., r_A$. Let us use lemma 3 and statement 15 of [2]. Let us rewrite these statements with respect to the case with two variables.

**Lemma 2.** Let $0_{(2)}$ be the point of local minimum of the polynomial $p(x, y)$, $p(0,0) = 0$, $p'(0,0) = 0_{(2)}$. Then $\forall A \in \mathbb{N}^2$ for the decomposition (4.1), (4.2) the following holds: $\forall i \in \{1, 2, ..., r_A\}$, $\forall (x, y) \in H_{\varphi_1^A} \cap ... \cap H_{\varphi_{i-1}^A}$ $\varphi_i^A(x, y) \geq 0$ (in particular for $i = 1$ we have $\forall (x, y) \in \mathbb{R}^2$ $\varphi_1^A(x, y) \geq 0$).

**Theorem 2.** Let $p(x, y)$ be a polynomial, $p(0,0) = 0$, $p'(0,0) = 0_{(2)}$, and $\forall A \in \mathbb{N}^2$ all main $A$-quasi-homogenious forms of the polynomial $p(x, y)$ are non-negative and nondegenerate in the weak sense (i.e., the necessary conditions of local minimality from Lemma 1 are satisfied). Let it be true for any $A \in \mathbb{N}_0^2$ for the decomposition (4.1), (4.2) that for some $j \in \{2, ..., r_A\}$ $0_{(2)}$ is a point of local minimum of the polynomials $\varphi_1^A(x, y)$, $\varphi_1^A(x, y) + \varphi_2^A(x, y)$,

..., $\varphi_1^A(x,y)+...+\varphi_{j-1}^A(x,y)$. In addition, either $H_{\varphi_1^A}\cap...\cap H_{\varphi_{j-1}^A}\cap[\mathbb{R}\setminus\{0\}]^2=\varnothing$ or

$\forall(x,y)\in H_{\varphi_1^A}\cap...\cap H_{\varphi_{j-1}^A}\cap[\mathbb{R}\setminus\{0\}]^2$ $\varphi_j^A(x,y)>0$. Then $0_{(2)}$ is the point of local minimum of $p(x,y)$.

As it follows from the above statements, in the case of $A_p\neq\varnothing$, other terms of the decomposition (4.1), (4.2) (besides $\varphi_1^A(x,y)$) and in particular $\varphi_2^A(x,y)$ may be needed for further investigations. Since $\varphi_2^A(x,y)$ is also a $A$-quasi-homogenious form, similarly to $\varphi_1^A(x,y)$ this form can be represented as

$$\varphi_2^A(x,y)=\sum_{i=1}^{v}b_i x^{\chi_i}y^{\eta_i}, \tag{4.3}$$

$$b_i\neq 0,\chi_i,\eta_i\in\mathbb{N}\cup\{0\}, A_1\chi_i+A_2\eta_i=B_2^A, i=1,...,v, v\in\mathbb{N},$$
$$\chi_1>\chi_2>...>\chi_v\geq 0. \tag{4.4}$$

The points $(\chi_i,\eta_i)$ are on the same line whose coordinates $(\chi,\eta)$ satisfy the equation $A_1\chi+A_2\eta=B_2^A$ or $\dfrac{\alpha-\chi_1}{-A_2}=\dfrac{\beta-\eta_1}{A_1}$ with a direction vector $e=(-A_2,A_1)$, and hence $(\chi_i,\eta_i)=(\chi_1,\eta_1)+\mu_i e$, $i=1,...,v$, where $\mu_1=0$, $\mu_i>0$, $i=2,...,v$. Using the fact that $\text{GCD}(|e_1|,|e_2|)=\text{GCD}(A_1,A_2)=1$, $(\chi_i-\chi_1,\eta_i-\eta_1)=\mu_i e$, we obtain that $\mu_i\in\mathbb{N}$, $i=2,...,v$ ( $\mu_1=0$ ). Thus, at $x\neq 0$ we have:

$$\varphi_2^A(x,y)=\sum_{i=1}^{v}b_i x^{\chi_i}y^{\eta_i}=\sum_{i=1}^{v}b_i x^{\chi_1+\mu_i e_1}y^{\eta_1+\mu_i e_2}=x^{\chi_1}y^{\eta_1}\sum_{i=1}^{v}b_i\left(x^{e_1}y^{e_2}\right)^{\mu_i}=$$
$$=x^{\chi_1}y^{\eta_1}\sum_{i=1}^{v}b_i u^{\mu_i}=x^{\chi_1}y^{\eta_1}g_2^A(u),\ u=x^{e_1}y^{e_2},\ g_2^A(u)=\sum_{i=1}^{v}b_i u^{\mu_i}, \tag{4.5}$$

where $g_2^A(u)$ is the characteristic polynomial of the quasi-homogenious form $\varphi_2^A(x,y)$, and $b_1 x^{\chi_1}y^{\eta_1}$ is its main term.

Consider the following two conditions:

(C1)$_A$ $\exists(x_0,y_0)\in H_{\varphi_1^A}:\ \varphi_2^A(x_0,y_0)<0$;

(C2)$_A$ $\forall(x,y)\in H_{\varphi_1^A}\ \left[x\neq 0, y\neq 0\Rightarrow\varphi_2^A(x,y)>0\right]$.

A consequence of Lemma 2 and Theorem 2 is

**Theorem 3.** Let $p(x,y)$ be a polynomial, $p(0,0)=0$, $p'(0,0)=0_{(2)}$, $\dim\text{Co}\,N_p=2$, and condition (3.1) is satisfied (i.e., all main quasi-homogenious forms from groups 1 and 2 are non-negative and nondegenerate in the weak sense)). Then, if $A_p=\varnothing$, then $0_{(2)}$ is the point of local minimum of the polynomial $p(x,y)$. If $A_p\neq\varnothing$, then the following cases are possible:

1) If $\exists A\in A_p$, that the condition (C1)$_A$ is satisfied, then $0_{(2)}$ is not a point of local minimum of the polynomial $p(x,y)$;

2) If for any $A\in A_p$, (C1)$_A$ is not satisfied, but the condition (C2)$_A$ is true, then $0_{(2)}$ is the point of local minimum of the polynomial $p(x,y)$.

Let us now give statements that reduce the verification of conditions (C1)$_A$, (C2)$_A$ to the solution of very simple computational problems. Let all main quasi-homogenious forms of the polynomial $p(x, y)$ be non-negative and $A_p \neq \emptyset$. Consider for $A \in A_p, u_0 \in U_p(A)$ the condition $(C3)_{A,u_0}$, which is satisfied if and only if the following system not joint:

$$\begin{cases} x \neq 0, y \neq 0, \\ x^{e_1} y^{e_2} = u_0, \\ \varphi_2^A(x, y) = x^{\chi_1} y^{\eta_1} g_2^A(u_0) < 0. \end{cases} \qquad (4.6)$$

**Assertion 3.** Suppose we are under the conditions of Theorem 3 and $A \in A_p \neq \emptyset$. Then it is necessary and sufficient for at least one of the three conditions to be true for condition $(C1)_A$ to hold (in the case of which $0_{(2)}$ is not a point of local minimum of the polynomial $p(x, y)$):

$$\alpha_s > 0, \; \chi_v = 0, \; (b_v < 0) \vee (\eta_v \notin 2\mathbb{N}), \qquad (4.7)$$

$$\beta_1 > 0, \; \eta_1 = 0, \; (b_1 < 0) \vee (\chi_1 \notin 2\mathbb{N}), \qquad (4.8)$$

$\exists u_0 \in U_p(A)$ such that $(У3)_{A,u_0}$ is not true
(i.e. the system (4.6) is joint). $\qquad (4.9)$

**Proof. Sufficiency.** In case (4.7) (case (4.8) is treated analogously), if $b_v < 0$, then at $x = 0, y = 1$ we have: $\varphi_1^A(0,1) = 0$, $\varphi_2^A(0,1) = b_v < 0$, and, if $b_v > 0, \eta_v \notin 2\mathbb{N}$ ($\eta_v \neq 0$, since $\chi_v = 0$), then at $x = 0, y = -1$ we have: $\varphi_1^A(0,-1) = 0$, $\varphi_2^A(0,-1) = -b_v < 0$. The case (4.9) is obvious (for $x, y$, satisfying (4.6), we have $(x, y) \in H_{\varphi_1^A}$, $\varphi_2^A(x, y) < 0$).

**Necessity.** Let $\exists (x_0, y_0) \in H_{\varphi_1^A} : \varphi_2^A(x_0, y_0) < 0$. Let us show that then at least one of the conditions (4.7) - (4.9) is true. Let us first consider the case when $x_0 \neq 0, y_0 \neq 0$. Then for $u_0 = x_0^{e_1} y_0^{e_2}$ it holds that: $\varphi_1^A(x_0, y_0) = x_0^{\alpha_1} y_0^{\beta_1} g_1^A(u_0)$, $\varphi_2^A(x_0, y_0) = x_0^{\chi_1} y_0^{\eta_1} g_2^A(u_0)$, and, since $(x_0, y_0) \in H_{\varphi_1^A}$, then $g_1^A(u_0) = 0$, and hence, $u_0 \in U_p(A)$. But then $x = x_0, y = y_0, u_0$ are solutions of the system (4.6), i.e., (4.9) is true. Let $x_0 = 0$. Then from the condition $\varphi_2^A(x_0, y_0) < 0$ (taking into account $B_2^A > 0$) we obtain that $y_0 \neq 0$. In this case $\alpha_s = 0$ cannot be fulfilled, since if $\alpha_s = 0$, then $\varphi_1^A(0, y) = a_s y^{\beta_s}$, whence from the condition $(x_0, y_0) = (0, y_0) \in H_{\varphi_1^A}$ we obtain $\varphi_1^A(0, y_0) = a_s y_0^{\beta_s} = 0$, and this contradicts the fact that $y_0 \neq 0$. Thus, $\alpha_s > 0$. Note now that $\chi_v > 0$ cannot be fulfilled, since then $\varphi_2^A(0, y_0) = b_v \cdot 0^{\chi_v} \cdot y_0^{\eta_v} = 0$. Thus, $\chi_v = 0$. Further, it cannot be simultaneously $b_v > 0$ and $\eta_v \in 2\mathbb{N}$, since then $\varphi_2^A(x_0, y_0) = \varphi_2^A(0, y_0) = b_v y_0^{\eta_v} > 0$. It is proved in a very similar way that in the case of $y_0 = 0$, (4.8) holds.

**Example 4.** 1) The case (4.7) corresponds to the example:
$$p(x, y) = \varphi_1^A(x, y) + \varphi_2^A(x, y) = x^2(x - y)^2 + 2y^5.$$
Here $A = (1, 1) \in \mathbb{N}^2$, $\varphi_1^A(x, y) = x^2(x - y)^2$, $\varphi_2^A(x, y) = 2y^5$, $\varphi_1^A(0, -1) = 0$, $\varphi_2^A(0, -1) = -2 < 0$.

2) The case (4.8) corresponds to the example:
$$p(x, y) = \varphi_1^A(x, y) + \varphi_2^A(x, y) = y^2(x - y)^2 + 2x^5.$$
Here $A = (1, 1) \in \mathbb{N}^2$, $\varphi_1^A(x, y) = y^2(x - y)^2$, $\varphi_2^A(x, y) = 2x^5$, $\varphi_1^A(-1, 0) = 0$, $\varphi_2^A(-1, 0) = -2 < 0$.

**Assertion 4.** Suppose we are under the conditions of Theorem 3, $A \in A_p \neq \emptyset$, and the condition $(C2)_A$ is satisfied. Then $\forall u_0 \in U_p(A)$ the condition $(C3)_{A,u_0}$ is satisfied.

**Proof.** Suppose that for some $u_0 \in U_p(A)$ the condition $(C3)_{A,u_0}$ is not satisfied, i.e. the system (4.6) is joint. Then it is true for $x, y$ satisfying this system:
$$\varphi_1^A(x,y) = x^{\alpha_1} y^{\beta_1} g_1^A(u_0) = 0, \ \varphi_2^A(x,y) = x^{\chi_1} y^{\eta_1} g_2^A(u_0) < 0, \ x \neq 0, y \neq 0,$$
and this contradicts the condition $(C2)_A$.

**Assertion 5.** Suppose we are under the conditions of Theorem 3, $A \in A_p \neq \emptyset$, and for any $u_0 \in U_p(A)$ the conditions $(C3)_{A,u_0}$, $g_2^A(u_0) \neq 0$ are satisfied. Then the condition $(C2)_A$ is true.

**Proof.** If condition $(C2)_A$ is not satisfied, then there are $(x,y) \in H_{\varphi_1^A}$ such that $x \neq 0, \ y \neq 0, \ \varphi_2^A(x,y) \leq 0$. But then $\varphi_1^A(x,y) = x^{\alpha_1} y^{\beta_1} g_1^A(u_0) = 0$, where $x^{e_1} y^{e_2} = u_0$, whence $g_1^A(u_0) = 0$, and hence $u_0 \in U_p(A)$. In this case, $\varphi_2^A(x,y) = x^{\chi_1} y^{\eta_1} g_2^A(u_0) \neq 0$, and hence $\varphi_2^A(x,y) = x^{\chi_1} y^{\eta_1} g_2^A(u_0) < 0$, i.e. the system (4.6) turns out to be joint, which contradicts the initial assumption that the condition $(C3)_{A,u_0}$ is true.

From statements 4, 5 we obtain that it is true that

**Corollary 1.** Suppose that we are in the conditions of Theorem 3, $A \in A_p \neq \emptyset$, $\forall u_0 \in U_p(A) \ g_2^A(u_0) \neq 0$. Then it is necessary and sufficient for condition $(C2)_A$ to hold for any $u_0 \in U_p(A)$ that condition $(C3)_{A,u_0}$ is satisfied (i.e., the system (4.6) is not joint).

Thus, checking for some $A \in A_p \neq \emptyset$ the conditions $(C1)_A$, $(C2)_A$ reduces to computing $g_2^A(u_0)$, $u_0 \in U_p(A)$, and checking (in the case of $g_2^A(u_0) \neq 0$) the jointness of the system (4.6).

The investigation of the jointness (or incompatibility) of the system (4.6) in the case $g_2^A(u_0) \neq 0$ is not difficult. Since $\text{GCD}(|e_1|, |e_2|) = 1$, at least one of the numbers among $e_1, e_2$ is odd. Suppose, for example, that the number $e_1$ is odd. Then for any value of $y \in \mathbb{R}$ from the second equality in (4.6) we uniquely determine the value of $x$ according to the formula $x = (u_0 y^{-e_2})^{1/e_1}$, substituting it into the third condition from (4.6), we obtain the inequality $(u_0 y^{-e_2})^{\chi_1/e_1} y^{\eta_1} g_2^A(u_0) < 0$ or

$$\left[ u_0^{\chi_1/e_1} g_2^A(u_0) \right] y^{\frac{e_1 \eta_1 - e_2 \chi_1}{e_1}} < 0. \tag{4.10}$$

If the integer $e_1 \eta_1 - e_2 \chi_1$ is even, then condition (4.10) is equivalent to the easily verifiable condition

$$u_0^{\chi_1/e_1} g_2^A(u_0) < 0. \tag{4.11}$$

If it is fulfilled, we make sure that the system (4.6) is joint. Otherwise, $u_0^{\chi_1/e_1} g_2^A(u_0) > 0$ is satisfied (recall that $u_0 \neq 0$) and the system (4.6) is incompatible. If the number $e_1 \eta_1 - e_2 \chi_1$ is odd, then at any value of $u_0^{\chi_1/e_1} g_2^A(u_0)$, which is different from 0 in the considered case, it is possible to choose the value of $y$ (for example, to choose it from the numbers 1, -1) so that the condition (4.10) holds, i.e., in this case the system (4.6) is joint (after choosing $y$, we assume $x = (u_0 y^{-e_2})^{1/e_1}$, so that the equality $x^{e_1} y^{e_2} = u_0$ holds).

Now let $e_2$ be odd. Then for any value of $x \in \mathbb{R}$ from the second equality in (4.6) we uniquely determine the value of $y$ according to the formula $y = (u_0 x^{-e_1})^{1/e_2}$, substituting it into the third condition in (4.6), we obtain the inequality $x^{\chi_1} (u_0 x^{-e_1})^{\eta_1/e_2} g_2^A(u_0) < 0$ or

$$\left[ u_0^{\eta_1/e_2} g_2^A(u_0) \right] x^{\frac{e_2 \chi_1 - e_1 \eta_1}{e_2}} < 0. \tag{4.12}$$

If the integer $e_2 \chi_1 - e_1 \eta_1$ (or, what is the same, $e_1 \eta_1 - e_2 \chi_1$; see the previous case) is even, then condition (4.12) is equivalent to the easily checked condition

$$u_0^{\eta_1/e_2} g_2^A(u_0) < 0, \tag{4.13}$$

in the case of which we make sure that the system (4.6) is joint. Otherwise, $u_0^{\eta_1/e_2} g_2^A(u_0) > 0$ and the system (4.6) is incompatible. If $e_2 \chi_1 - e_1 \eta_1$ is odd, then for any value of $u_0^{\eta_1/e_2} g_2^A(u_0)$, which in this case is different from 0, we can choose a value of $x$ (for example, choose it from the numbers 1, -1) so that the condition (4.12) holds, i.e., in this case the system (4.6) is joint (after choosing $x$, we assume $y = (u_0 x^{-e_1})^{1/e_2}$ so that the equality $x^{e_1} y^{e_2} = u_0$ holds). Acting in this way, for each $u_0 \in U_p(A)$, at which $g_2^A(u_0) \neq 0$, we uniquely determine whether the system (4.6) is joint or not.

On the basis of these results, we can already describe a rather simple algorithm for checking for a polynomial (or a power series; see Remark 8 below) $p(x, y) \neq 0$, where $p(0,0) = 0$, $p'(0,0) = 0_{(2)}$ (i.e., $0_{(2)}$ is a stationary point), whether $0_{(2)}$ is a point of local minimum of this polynomial. We will consider the nontrivial case when $\dim \mathrm{CoN}_p = 2$ and the single terms of the polynomial $p(x, y)$, corresponding to points from $\Omega(N_p)$ are non-negative (and hence nondegenerate in a weak sense quasi-homogeneous forms, i.e., (3.1) is satisfied.

**Algorithm 1**

**Step 1.** Identify all main quasi-homogeneous forms of the polynomial $p(x, y)$ from group 3 (see the definition of group 3 in Section 2 after lemma 1). If the set of such forms is empty, then by virtue of Theorem 1, the point $0_{(2)}$ is the point of local minimum of the polynomial $p(x, y)$ under consideration, and this is the end of the algorithm. Otherwise, suppose $\tilde{A}_p = \emptyset$, where $\tilde{A}_p$ is the set of all $A \in \mathbb{N}_0^2$ for which finer studies are required, and proceed to step 2.

**Step 2.** Choose any next main quasi-homogeneous form from group 3. If all of them have already been considered, then in the case $\tilde{A}_p = \emptyset$ (everywhere in the algorithm we proceed to step 2 in one of three cases: either in the case when the next quasi-homogeneous form $\varphi_1^A(x, y)$ from group 3 is non-negative and nondegenerate in the weak sense, i.e., when $A \notin A_p$, or if $A \in A_p$, and the condition $(C1)_A$ is not satisfied, but $(C2)_A$ is true, i.e., the condition 2) of Theorem 3 is not violated, or in the case of $A \in \tilde{A}_p$) the point $0_{(2)}$ is the point of local minimum of the polynomial $p(x, y)$ under consideration, and if $\tilde{A}_p \neq \emptyset$, then more subtle investigations are required (see below the modification of Algorithm 1, Remark 7, and Section 5), and this is the end of the algorithm's work. Otherwise, we define by any two terms of the chosen form a vector $A \in \mathbb{N}_0^2$ such that this form is the main $A$-quasi-homogeneous form of the polynomial $p(x,y)$. Since $\dim \text{Co} \, N_p = 2$, $p(x,y)$ contains at least two terms in the decomposition (4.1), (4.2), the first of which is $\varphi_1^A(x,y)$ and corresponds to the next chosen form from group 3.

**Step 3.** By virtue of (3.1), we obtain that condition 1) of statement 1 holds for $\varphi_1^A(x,y)$ (using statement 10 and Remark 9 of [2], it is easy to show that $(\alpha_1, \beta_1)$, $(\alpha_s, \beta_s) \in \Omega(N_p) = \Psi(N_p) \cap P(N_p)$). We find the characteristic polynomial $g_1^A(u)$ for the quasi-homogeneous form $\varphi_1^A(x,y)$ according to formula (2.3), and also the set $U_p(A)$. If $U_p(A) = \varnothing$, then the quasi-homogeneous form $\varphi_1^A(x,y)$ is nonnegative and nondegenerate in the weak sense by virtue of statement 2. Then $A \notin A_p$ and proceed to step 2. Otherwise $U_p(A) \neq \varnothing$ and proceed to step 4.

**Step 4.** Find multiplicities of roots from $U_p(A)$, by which we check the non-negativity of the polynomial $g_1^A(u)$. If it is not satisfied, then by virtue of statement 1 the quasi-homogeneous form $\varphi_1^A(x,y)$ is not non-negative, and consequently, by virtue of lemma 1, the point $0_{(2)}$ is not a point of local minimum of the polynomial $p(x,y)$, and this is the end of the algorithm. Otherwise $A \in A_p$, and we proceed to step 5.

**Step 5.** Check the fulfillment of the condition $(C1)_A$. To do this, we first determine the validity of (4.7), (4.8). If at least one of these conditions is satisfied, then by virtue of statement 3, the point $0_{(2)}$ is not the point of local minimum of the polynomial $p(x,y)$ and this is the end of the algorithm. Otherwise, we proceed to step 6.

**Step 6.** Find the characteristic polynomial $g_2^A(u)$ of the quasi-homogeneous form $\varphi_2^A(x,y)$ according to formula (4.5). If $\forall u_0 \in U_p(A) \; g_2^A(u_0) = 0$, we proceed to step 7. Otherwise, we consider two cases:

*Case 1*. Let the integer $e_1 \eta_1 - e_2 \chi_1$ be odd. Then the system (4.6) is joint for every $u_0 \in U_p(A)$, for which $g_2^A(u_0) \neq 0$, and according to statement 3, the point $0_{(2)}$ is not the point of the local minimum of the polynomial $p(x,y)$. This completes the work of the algorithm.

*Case 2*. Let $e_1 \eta_1 - e_2 \chi_1$ be an even integer. Then for each $u_0 \in U_p(A)$, for which $g_2^A(u_0) \neq 0$, we check whether condition (4.11) is satisfied if the integer $e_1$ is odd, or condition (4.13) is satisfied if the integer $e_2$ is odd (if both are odd, these conditions are equivalent). In the first case of simultaneous fulfillment of (4.11) and oddness $e_1$ or simultaneous fulfillment of (4.13) and oddness $e_2$ the system (4.6) is joint, then by virtue of statement 3 the point $0_{(2)}$ is not the point of the local minimum of the polynomial $p(x,y)$, and this is the end of the algorithm. If $\forall u_0 \in U_p(A) \; g_2^A(u_0) \neq 0$, and the system (4.6) is not joint (which is uniquely determined by checking one of the conditions (4.11) or (4.13)), then, by virtue of Corollary 1, $(C2)_A$ is satisfied, and by virtue of Assertion 3, $(C1)_A$ is not satisfied, i.e., condition 2) of Theorem 3 is not violated for the considered $A \in A_p$. In this case we consider the next main quasi-homogeneous form from group 3, i.e. we proceed to step 2. If we find a number $u_0 \in U_p(A)$ such that $g_2^A(u_0) = 0$, then we proceed to step 7.

**Step 7.** At this step we find ourselves in the case when $A \in A_p$, $g_2^A(u_0) = 0$ for some (at least one) $u_0 \in U_p(A)$, and the system (4.6) is incompatible for those $u_0 \in U_p(A)$ for which $g_2^A(u_0) \neq 0$. Then even more subtle investigations will be needed for a given $A \in A_p$. Therefore, we assign $\tilde{A}_p := \tilde{A}_p \cup \{A\}$. In any case, we can continue the algorithm by investigating the next main quasi-homogeneous form from group 3 (since this investigation may lead to a situation showing that $0_{(2)}$ is not a point of local minimum of the polynomial $p(x, y)$), i.e., we proceed to step 2.

Applying Algorithm 1 to a polynomial $p(x, y)$ of the considered kind, we either determine whether $0_{(2)}$ is a point of local minimum of this polynomial or not, or we select a non-empty set $\tilde{A}_p$ of all $A \in A_p$, for which more subtle investigations are required. Let us describe some of them.

These investigations can, in particular, be based on the application of Theorem 2 for cases with $j \geq 3$, i.e., using other terms in the decomposition (4.1), (4.2), in particular, the quasi-homogeneous form $\varphi_3^A(x, y)$. We will describe a modification of Algorithm 1 that uses $\varphi_3^A(x, y)$. Similarly, modifications using $\varphi_4^A(x, y)$ etc. can be described.

**Remark 6.** Note that there is a possible case when $p(x, y) = \varphi_1^A(x, y) + \varphi_2^A(x, y)$, i.e., $p(x, y)$ is the sum of the first two quasi-homogeneous forms. This case will be considered in Section 6 and a procedure will be derived for it which allows us to uniquely answer the question whether $0_{(2)}$ is a point of local minimum of this polynomial or not. In this connection, let us proceed to the case when, $\forall A \in \tilde{A}_p$ $\varphi_3^A(x, y) \not\equiv 0$.

By analogy with the condition $(C2)_A$, we introduce the condition (see also Remark 7 below)

$(\tilde{C}2)_A$ $\forall (x, y) \in H_{\varphi_1^A} \cap H_{\varphi_2^A}$ $\left[ x \neq 0, y \neq 0 \Rightarrow \varphi_3^A(x, y) > 0 \right]$; $\forall (x, y) \in \mathbb{R}^2$ $\varphi_2^A(x, y) \geq 0$.

Note that the second condition in $(\tilde{C}2)_A$ is checked in accordance with statement 1, by analogy with $\varphi_1^A(x, y)$. This condition can be replaced by a weaker one (see Remark 7 below).

A corollary of Lemma 2 and Theorem 2 (at $j = 3$) is

**Theorem 4.** Suppose we are under the conditions of Theorem 3. Then condition 2) in Theorem 3 can be replaced by the condition

2′) If $\forall A \in A_p$ satisfies the condition $(C2)_A$ or $(\tilde{C}2)_A$, then $0_{(2)}$ is the point of local minimum of the polynomial $p(x, y)$.

Let us denote by $\tilde{U}_p(A) = \{ u \in U_p(A) | g_2^A(u) = 0 \} = \{ u \in \mathbb{R} | g_1^A(u) = 0, g_2^A(u) = 0 \}$. By analogy with condition $(C3)_{A,u_0}$, consider for $A \in A_p \neq \emptyset$, $u_0 \in \tilde{U}_p(A)$ the condition $(\tilde{C}3)_{A,u_0}$ which is satisfied if and only if the system is not joint:

$$\begin{cases} x \neq 0, y \neq 0, \\ x^{e_1} y^{e_2} = u_0, \\ \varphi_3^A(x, y) = x^{\omega_1} y^{\xi_1} g_3^A(u_0) < 0, \end{cases} \quad (4.14)$$

obtained by analogy with the system (4.6) using the quasi-homogeneous form $\varphi_3^A(x, y)$, and also $g_3^A(u)$ that is its characteristic polynomial of $\varphi_3^A(x, y)$ (here $\omega_1$, $\xi_1$ are the degrees of the variables $x$, $y$ in the main term of the form $\varphi_3^A(x, y)$).

Quite analogously to assertions 4, 5, it is not difficult to show that are true:

**Assertion 6.** Suppose we are under the conditions of Theorem 3, $A \in A_p \neq \varnothing$, and the condition $(\tilde{C}2)_A$ is satisfied. Then $\forall u_0 \in \tilde{U}_p(A)$ the condition $(\tilde{C}3)_{A,u_0}$ is true.

**Assertion 7.** Let $\forall u_0 \in \tilde{U}_p(A)$ the following conditions be satisfied: $(\tilde{C}3)_{A,u_0}$, $g_2^A(u_0) \neq 0$. Then the condition $(\tilde{C}2)_A$ is true.

Using Theorem 4, as well as assertions 6, 7, we can modify Algorithm 1 by replacing step 7 with a new step $7'$ (this results in an additional step 8).

**Step $7'$.** In this step we find ourselves in the case where $A \in A_p$, $g_2^A(u_0) = 0$ is satisfied for some (at least one) $u_0 \in U_p(A)$, and yet the system (4.6) is incompatible for those $u_0 \in U_p(A)$ for which $g_2^A(u_0) \neq 0$. In this case we will need to consider a quasi-homogeneous form $\varphi_3^A(x, y)$. If $\varphi_3^A(x, y) \equiv 0$, then, in accordance with Remark 6, it is possible, using the procedure specified therein, to determine unambiguously whether $0_{(2)}$ is a point of local minimum of the polynomial $p(x, y)$ or not. Otherwise, we investigate for non-negativity of $\varphi_2^A(x, y)$ similarly to the investigation for non-negativity of $\varphi_1^A(x, y)$, described in steps 3-5 of Algorithm 1 (see also Remark 7 below). If it is not satisfied, even more subtle investigations will be needed for a given $A \in A_p$ (see Section 5). Therefore, we assign $\tilde{A}_p := \tilde{A}_p \cup \{A\}$ and proceed to step 2. Otherwise, the second condition in $(\tilde{C}2)_A$ is satisfied and we check the fulfillment of the first one. We proceed similarly to step 6, applying all the steps described in this step to the quasi-homogeneous form $\varphi_3^A(x, y)$ (instead of $\varphi_2^A(x, y)$) with all the resulting reassignments. If for all $u_0 \in \tilde{U}_p(A)$ $g_3^A(u_0) \neq 0$ and the system (4.14) is not joint, then, by virtue of statement 7, the condition $(\tilde{C}2)_A$ is true, i.e., the condition $2'$) of Theorem 4 is not violated for the $A \in \mathbb{N}_0^2$ under consideration. In this case we consider the next main quasi-homogeneous form from group 3, i.e. we proceed to step 2.

**Step 8.** At this step we find ourselves in the case when $A \in A_p$, and $g_3^A(u_0) = 0$ for some (at least one) $u_0 \in \tilde{U}_p(A)$, but the system (4.14) is incompatible for those $u_0 \in \tilde{U}_p(A)$ for which $g_3^A(u_0) \neq 0$. In this case, we may need to consider a quasi-homogeneous form $\varphi_4^A(x, y)$ and further modify the algorithm up to the exhaustion of all quasi-homogeneous forms in the decomposition (4.1), (4.2). Note, however, that further modification of the algorithm is possible, for example, if the non-negativity condition is satisfied (see also Remark 7 concerning the relaxation of this condition). Otherwise, we assign $\tilde{A}_p := \tilde{A}_p \cup \{A\}$ and proceed to step 2.

**Remark 7.** In Theorem 2 we use the condition that for some $j \in \{2, \ldots, r_A\}$ the point $0_{(2)}$ is the point of local minimum of the polynomials $\varphi_1^A(x), \varphi_1^A(x) + \varphi_2^A(x), \ldots, \varphi_1^A(x) + \cdots + \varphi_{j-1}^A(x)$. This condition is known to be satisfied if the easily verifiable non-negativity condition $\varphi_1^A(x), \ldots, \varphi_{j-1}^A(x)$ is true (thus, by virtue of $A \in \mathbb{N}^2$, the non-negativity condition of a quasi-homogeneous form $\varphi_1^A(x)$ is equivalent to the condition: $0_{(2)}$ is the point of local minimum of $\varphi_1^A(x)$). Section 6 describes a procedure that allows us to uniquely answer the question whether $0_{(2)}$ is a point of local minimum of the polynomial $\varphi_1^A(x)+\varphi_2^A(x)$ or not. In this connection, the second condition in $(\tilde{C}2)_A$ (i.e., the non-negativity condition $\varphi_2^A(x,y)$) can be replaced by a check of the local minimum of the function $\varphi_1^A(x)+\varphi_2^A(x)$ at the point $0_{(2)}$. But then in the modification of Algorithm 1 we can replace the check of the non-negativity condition $\varphi_2^A(x,y)$ by the check of the weaker condition of the local minimum of the function $\varphi_1^A(x)+\varphi_2^A(x)$ at the point $0_{(2)}$.

Accordingly, in further modification of Algorithm 1, we can continue to replace the non-negativity conditions of polynomials $\varphi_3^A(x,y), \varphi_4^A(x,y)$, etc. by weaker conditions: $0_{(2)}$ is the point of local minimum of polynomials $\varphi_1^A(x,y)+\varphi_2^A(x,y)+\varphi_3^A(x,y)$, $\varphi_1^A(x,y)+\varphi_2^A(x,y)+\varphi_3^A(x,y)+\varphi_4^A(x,y)$, etc.

**Remark 8.** All the statements used in Algorithm 1 and its modification remain valid for the power series $p(x,y)$, absolutely convergent in some neighborhood of the point $0_{(2)}$ (see [5]), and hence these algorithms can be applied to $p(x,y)$ in this case as well. In this case, we find the set of main $A$-quasi-homogeneous forms where $A \in \mathbb{N}_0^2$, and then the set of vectors of $\tilde{A}_p$ by $P(N_p)$, which is finite (see [5]). Thus, the quasi-homogeneous forms $\varphi_1^A(x)$, where $A \in A_p$, are found by the terms of the series $p(x,y)$ corresponding to the finite set $P(N_p)$. Further, $\varphi_2^A(x,y)$ may require terms of the power series $p(x,y)$, corresponding to a finite set $P(N_p \setminus N_{\varphi_1^A})$, etc.

**Example 5.** Consider a polynomial $p(x,y,a)$ (see $N_p$, $\mathrm{Co}\, N_p$ in Fig. 1) for which the following holds at $A=(1,\ 2)$

$$p(x,y,a) = \varphi_1^A(x,y)+\varphi_2^A(x,y,a)+\varphi_3^A(x,y),$$

where $\varphi_1^A(x,y) = x^4 y^2 + 2x^2 y^3 + y^4 = y^2(x^2+y)^2$, $\varphi_2^A(x,y,a) = 3ax^6 y^2 + 3x^4 y^3$, $\varphi_3^A(x,y) = 0.01 x^8 y^3$ are polynomials which are $A$-quasi-homogeneous forms, $B_1^A = 8, B_2^A = 10$, $B_3^A = 14$ (see (4.1), (4.2)), $a$ is the parameter ($a \in \{0.99,\ 1,\ 1.01\}$). In this example $p(0,0,a)=0$, $p'(0,0,a)=0_{(2)}$, i.e. $0_{(2)}$ is a stationary point. Let us apply Algorithm 1 to $p(x,y,a)$. It is easy to see that for this polynomial at any of the given values of $a$ the condition (3.1) holds, i.e. all the main quasi-homogeneous forms of this polynomial from groups 1, 2 are non-negative, nondegenerate in the weak sense and the only main quasi-homogeneous form from group 3 is $\varphi_1^A(x,y)$. The vector $A=(1,\ 2) \in \mathbb{N}_0^2$ corresponding to it is uniquely determined by the members of this form (any two of them). At steps 3, 4 we investigate the non-negativity of $\varphi_1^A(x,y)$. Note that $x^4 y^2$ is its main term, $\varphi_1^A(x,y) = x^4 y^2 g_1^A(u)$, where $g_1^A(u) = (1+u)^2$ is its characteristic polynomial, $u = x^{-2} y = x^{e_1} y^{e_2}$, $e_1 = -A_2 = -2$, $e_2 = A_1 = 1$. Thus, the characteristic polynomial $g_1^A(u)$ has a single real root $u_0 = -1$ of multiplicity 2, i.e., $U_p(A) = \{-1\}$. In this case, the conditions of statement 1 are satisfied, i.e., the quasi-homogeneous form $\varphi_1^A(x, y)$ is nonnegative. Thus, $A_p = \{(1, 2)\}$. In step 5, we check the fulfillment of conditions (4.7), (4.8), which are obviously not satisfied. At step 6 we study the quasi-homogeneous form $\varphi_2^A(x, y, a) = x^{\chi_1} y^{\eta_1} g_2^A(u, a)$, where $b_1 x^{\chi_1} y^{\eta_1} = 3x^6 y^2$ is its main term and $g_2^A(u, a) = 3(a+u)$ is its characteristic polynomial. Note that the number $e_1 \eta_1 - e_2 \chi_1 = -2 - 10 = -12$ is even. Therefore, according to step 6, since the number $e_2 = A_1 = 1$ is odd, we check the fulfillment of condition (4.13), which is of the form:

$$u_0^{\eta_1/e_2} g_2^A(u_0, a) < 0 \Leftrightarrow g_2^A(-1, a) < 0 \Leftrightarrow 3(a-1) < 0.$$

In the case of $a = 0.99$, this condition is satisfied, and hence $0_{(2)}$ is the point of local minimum of the polynomial $p(x, y, 0.99)$. In the case of $a = 1.01$ it is not satisfied, i.e., the condition $(C3)_{A,u_0}$ is true. But then, by Corollary 1, $(C2)_A$ holds, and by Proposition 3, $(C1)_A$ does not hold. In this case, the uniqueness of $A = (1, 2) \in A_p$ implies the validity of condition 2) of Theorem 3, by virtue of which $0_{(2)}$ is the point of the local minimum of the polynomial $p(x, y, 1.01)$. In the case of $a = 1$, $g_2^A(u_0, 1) = 0$ is satisfied, i.e., this case requires more subtle investigations (in particular, those proposed in Section 5). It is not difficult to see that in the third case ($a = 1$), by substituting $x(t) = t$, $y(t) = -t^2$ the equality $p(x(t), y(t), 1) = -0.01 t^{14}$ is satisfied, i.e., $0_{(2)}$ is not a point of local minimum of the polynomial $p(x, y, 1)$. Note that in the third case the method of polynomial substitution, which is discussed in detail in Section 5, was actually applied.

## 5. METHOD OF SUBSTITUTION OF POLYNOMIALS WITH INDEFINITE COEFFICIENTS

From the consideration of Algorithm 1 and its modification, we can see that applying it to some polynomials leads to the case where the result of the work is a non-empty set $\tilde{A}_p$, showing that the question under investigation remains open and a more subtle investigation of the decompositions of the polynomial into $A$-quasi-homogeneous forms for $A \in \tilde{A}_p$ is required. This is because Algorithm 1 and its modification are based on the statements of Lemma 2 and Theorem 2, which give separately a necessary condition (Lemma 2) and separately a sufficient condition (Theorem 2) that are not combined into a single general condition of local minimum.

In this section, we consider a method that is based on a general necessary and sufficient condition for the local minimum of a polynomial $p(x, y)$. In this method, we consider polynomials with a variable $t \in \mathbb{R}$ of the form

$$x(t) = c_0 t^{\nu_1} + c_1 t^{\nu_1+1} + o(t^{\nu_1+1}), \ y(t) = d_0 t^{\nu_2} + d_1 t^{\nu_2+1} + o(t^{\nu_2+1}),$$
$$c_0, d_0 \neq 0, \ \nu_1, \nu_2 \in \mathbb{N}. \tag{5.1}$$

We will choose polynomials $x(t)$, $y(t)$, of the form (5.1) such that

$$p(x(t), y(t)) = g_0 t^\sigma + o(t^\sigma), \ \sigma \in \mathbb{N}, \ g_0 < 0. \tag{5.2}$$

In this case we will use the following assertion 8, which is a consequence of lemma 7 of [1].

**Assertion 8.** Let $p(x, y)$ be a polynomial, $p(x, y) \not\equiv 0$, $p(0,0) = 0$, $p'(0,0) = 0_{(2)}$, i.e., $0_{(2)}$ is a stationary point. For $0_{(2)}$ not to be a point of local minimum of the polynomial $p(x, y)$, it is necessary and sufficient that there exist polynomials of the form (5.1) such that (5.2) holds.

Note that for $A = (A_1, A_2) = (\nu_1, \nu_2)/\nu_0 \in \mathbb{N}_0^2$, where $\nu_0 = \mathrm{GCD}(\nu_1, \nu_2)$, for the decomposition (4.1), (4.2) it follows that

$$p(x(t), y(t)) = \varphi_1^A(x(t), y(t)) + \varphi_2^A(x(t), y(t)) + \ldots + \varphi_{r_A}^A(x(t), y(t)),$$

where $\varphi_1^A(x, y)$ is the main $A$-quasi-homogeneous form of the polynomial $p(x, y)$ satisfying (2.1) - (2.3), and hence,

$$p(x(t), y(t)) = \varphi_1^A(c_0, d_0) t^{B_1^A} + o(t^{B_1^A}). \tag{5.3}$$

The method of indeterminate coefficients applies in the nontrivial case when the main $A$-quasi-homogeneous forms of the polynomial $p(x, y)$ at $A \in \mathbb{N}_0^2$ are non-negative, i.e., in particular, the condition (3.1) is satisfied and thus $A_p = \varnothing$. If $A_p = \varnothing$, then, as shown in Section 2, the question of whether $0_{(2)}$ is a point of local minimum of the polynomial $p(x, y)$ is solved simply using Statement 2 and Theorem 1, i.e., in this case $0_{(2)}$ is a point of local minimum of the polynomial $p(x, y)$.

It follows from the non-negativity of $\varphi_1^A(x, y)$ that, by virtue of (5.3), condition (5.2) can be satisfied only when $\varphi_1^A(c_0, d_0) = 0$, and therefore only this case is of interest for what follows, in which, by virtue of (2.1) - (2.3), we have:

$$\begin{cases} \varphi_1^A(c_0, d_0) = c_0^{\alpha_1} d_0^{\beta_1} g_1^A(u_0) = 0, \ u_0 = c_0^{e_1} d_0^{e_2} \in U_p(A), \\ A = (A_1, A_2) = (\nu_1, \nu_2)/\nu_0 \in A_p, \ \nu_0 = \mathrm{GCD}(\nu_1, \nu_2), \ e_1 = -A_2, e_2 = A_1. \end{cases} \tag{5.4}$$

Thus, it was possible to narrow down the set of considered polynomials of the form (5.1) checked for fulfillment of condition (5.2). The set $A_p$ is finite, since the number of vectors in $A_p$ does not exceed the number of main quasi-homogeneous forms of the polynomial $p(x, y)$, belonging to group 3. Often $A_p$ consists of a single vector.

Consider a condition formulated with respect to $A \in A_p$, $u_0 \in U_p(A)$.

(C4)$_{A,u_0}$ There exist $c_0 \neq 0$, $d_0 \neq 0$ such that conditions (5.2), (5.4) are satisfied for polynomials of the form (5.1).

Let us also formulate a condition for a polynomial $p(x, y)$ such that $A_p \neq \varnothing$.

(C4) There exist a vector $A \in A_p$ and a number $u_0 \in U_p(A)$ such that the condition (C4)$_{A,u_0}$ holds.

The consequence of the above reasoning, as well as of statement 8 is

**Assertion 9.** Let $p(x, y)$ be a polynomial, $p(x, y) \not\equiv 0$, $p(0, 0) = 0$, $p'(0, 0) = 0_{(2)}$, i.e., $0_{(2)}$ is a stationary point. Let, further, all main quasi-homogeneous forms of the polynomial $p(x, y)$ of groups 1 and 2 be non-negative and nondegenerate in the weak sense, and thus $A_p \neq \varnothing$. Then for $0_{(2)}$ not to be a point of local minimum of the polynomial $p(x, y)$ it is necessary and sufficient that condition (C4) is satisfied.

For convenience in applying statement 9, for each $A \in A_p$ it is desirable to narrow down the set of polynomials of the form (5.1) that are checked in condition (C4)$_{A,u_0}$ for fulfillment of (5.4).

Let $A \in A_p$, $u_0 \in U_p(A)$, and we are going to check the fulfillment of condition $(C4)_{A,u_0}$. Note that condition $u_0 = c_0^{e_1} d_0^{e_2}$ is satisfied by an infinite set of pairs $c_0, d_0$. Let us narrow this set to a finite one. Note that by virtue of $\mathrm{GCD}(|e_1|, |e_2|) = 1$, at least one of the numbers among $e_1$, $e_2$ is odd. Thus, there are three possible cases, for each of which we consider the corresponding sets:

1) Let $e_1 = -A_2$ be even and $e_2 = A_1$ be odd. Let's assume $C_1 = \{(1, u_0^{1/e_2}), (-1, u_0^{1/e_2})\}$.

2) Let $e_1$ be odd and $e_2$ be even. Let's assume $C_2 = \{(u_0^{1/e_1}, 1), (u_0^{1/e_1}, -1)\}$.

3) Let $e_1$, $e_2$ be odd. Let's assume $C_3 = \{(1, u_0^{1/e_2}), (-1, -u_0^{1/e_2})\}$, $C_4 = \{(u_0^{1/e_1}, 1), (-u_0^{1/e_1}, -1)\}$.

Let us show that in any of these cases it is possible, when checking condition $(C4)_{A,u_0}$, to restrict ourselves to considering selected finite sets of pairs $(c_0, d_0)$.

Preliminarily we show that in each of these cases, if $u_0 = c_0^{e_1} d_0^{e_2}$ for some $c_0 \neq 0$, $d_0 \neq 0$, then we find (and easily compute) $\tau > 0$, such that the following is true for a pair of numbers $\tilde{c}_0 = c_0 \tau^{1/A_2}, \tilde{d}_0 = d_0 \tau^{1/A_1}$. In case 1) $\{(\tilde{c}_0, \tilde{d}_0), (-\tilde{c}_0, \tilde{d}_0)\} = C_1$. In case 2) $\{(\tilde{c}_0, \tilde{d}_0), (\tilde{c}_0, -\tilde{d}_0)\} = C_2$. In case 3) there are two values for a number $\tau > 0$ at one of which $\{(\tilde{c}_0, \tilde{d}_0), (-\tilde{c}_0, -\tilde{d}_0)\} = C_3$, and at the other $\{(\tilde{c}_0, \tilde{d}_0), (-\tilde{c}_0, -\tilde{d}_0)\} = C_4$.

Let us consider case 1) (the others are treated similarly). Indeed, in case 1) it is sufficient to put $\tau = c_0^{-A_2} = c_0^{e_1} = |c_0|^{e_1}$ (thus $|c_0| = \tau^{-1/A_2}$, $\tau = |c_0|^{-A_2}$). Then, taking into account the fact that, by virtue of $u_0 = c_0^{e_1} d_0^{e_2} = |c_0|^{e_1} d_0^{e_2}$, $(u_0)^{\frac{1}{A_1}} = (u_0)^{\frac{1}{e_2}} = d_0 |c_0|^{\frac{e_1}{e_2}} = d_0 |c_0|^{-\frac{A_2}{A_1}}$ is satisfied, we obtain that for $\tilde{c}_0 = \pm c_0 \tau^{1/A_2}, \tilde{d}_0 = d_0 \tau^{1/A_1}$ it is true: $|\tilde{c}_0| = |c_0| \tau^{1/A_2} = 1$,

$\tilde{d}_0 = d_0 \tau^{1/A_1} = d_0 |c_0|^{-\frac{A_2}{A_1}} = (u_0)^{\frac{1}{A_1}} = (u_0)^{\frac{1}{e_2}}$, i.e., $\{(\tilde{c}_0, \tilde{d}_0), (-\tilde{c}_0, \tilde{d}_0)\} = \{(1, u_0^{1/e_2}), (-1, u_0^{1/e_2})\} = C_1$ is satisfied.

Let us return to the problem of checking the condition $(C4)_{A,u_0}$ for some $A \in A_p$, $u_0 \in U_p(A)$. Let us show that this check is sufficient for a finite number of sets $c_0 \neq 0$, $d_0 \neq 0$, such that $u_0 = c_0^{e_1} d_0^{e_2} \in U_p(A)$.

Let us take advantage of the fact that for a vector $A \in A_p$ one of the cases listed above concerning the parity or oddness of its components is satisfied. Let, for example, case 1) (when $e_1 = -A_2$ is even and $e_2 = A_1$ is odd) be satisfied. The other cases are treated similarly. Suppose that for some $u_0 \in U_p(A)$ the condition $(C4)_{A,u_0}$ is satisfied. Then there exist $c_0 \neq 0$, $d_0 \neq 0$ such that conditions (5.2), (5.4) are true for polynomials of the form (5.1). Consider the substitution: $t = (\tau^{1/(v_0 A_1 A_2)}) \tilde{t}$, where $\tau > 0$. Then, taking into account (5.2), we obtain

$$p(x(t), y(t)) = p\left(x\left(\tau^{\frac{1}{v_0 A_1 A_2}} \tilde{t}\right), y\left(\tau^{\frac{1}{v_0 A_1 A_2}} \tilde{t}\right)\right) = \tilde{g}_0 \tilde{t}^\sigma + o(\tilde{t}^\sigma),$$

where $\tilde{g}_0 = g_0 \tau^{\frac{\sigma}{v_0 A_1 A_2}} < 0$. Moreover

$$x(t) = x\left(\tau^{\frac{1}{v_0 A_1 A_2}} \tilde{t}\right) = c_0 \tau^{\frac{v_1}{v_0 A_1 A_2}} \tilde{t}^{v_1} + c_1 \tau^{\frac{v_1+1}{v_0 A_1 A_2}} \tilde{t}^{v_1+1} + o(\tilde{t}^{v_1+1}),$$

$$y(t) = y\left(\tau^{\frac{1}{v_0 A_1 A_2}} \tilde{t}\right) = d_0 \tau^{\frac{v_2}{v_0 A_1 A_2}} \tilde{t}^{v_2} + d_1 \tau^{\frac{v_2+1}{v_0 A_1 A_2}} \tilde{t}^{v_2+1} + o(\tilde{t}^{v_2+1}).$$

Note that $\tilde{c}_0 = c_0 \tau^{v_1/(v_0 A_1 A_2)} = c_0 \tau^{1/A_2}$, $\tilde{d}_0 = d_0 \tau^{v_2/(v_0 A_1 A_2)} = d_0 \tau^{1/A_1}$, and hence, when $\tau = c_0^{-A_2} = c_0^{e_1} = |c_0|^{e_1}$ is chosen (as was shown in case 1)), $(\tilde{c}_0, \tilde{d}_0) \in C_1 = \{(1, u_0^{1/e_2}), (-1, u_0^{1/e_2})\}$ is satisfied, i.e., $\tilde{d}_0 = u_0^{1/e_2}$ is uniquely chosen, but $\tilde{c}_0 \in \{1, -1\}$. Thus, it is shown that when checking for any $A \in A_p$, $u_0 \in U_p(A)$ the condition (C4)$_{A,u_0}$ in case 1), an additional condition $(c_0, d_0) \in C_1$ can be imposed on polynomials of the form (5.1) satisfying (5.2), which allows us to restrict ourselves to considering only two of cases out of an infinite number of cases where the condition (C4)$_{A,u_0}$ is satisfied.

The situation is similar for the other two cases. Thus, instead of the condition (C4)$_{A,u_0}$ we can consider the equivalent condition

(C̃4)$_{A,u_0}$ There exist $c_0 \neq 0$, $d_0 \neq 0$, such that for polynomials of the form (5.1), conditions (5.2), (5.4) are satisfied and in case 1) $(c_0, d_0) \in C_1$, in case 2) $(c_0, d_0) \in C_2$, and in case 3) $(c_0, d_0) \in C_3$ or $(c_0, d_0) \in C_4$ (choose any of these conditions).

Let us now consider the question of choosing a vector $(v_1, v_2) \in \mathbb{N}^2$ for polynomials of the form (5.1) while checking the condition (C4)$_{A,u_0}$ or (C̃4)$_{A,u_0}$. One of the constraints is condition (5.4), under which we must find a number $v \in \mathbb{N}$ such that $(v_1, v_2) = vA$. This raises the question: is it possible to dispense with the case of $v = 1$? The following example shows that when applying the method of substitution of polynomials with undefined coefficients, it is not always possible to limit oneself to the value of $v = 1$.

**Example 6.** Let $p(x, y) = (x-y)^6 - (x-y)^2 x^5 + x^8$. Then $p(x, y) = \varphi_1^A(x, y) + \varphi_2^A(x, y) + \varphi_3^A(x, y)$ is satisfied for $A = (A_1, A_2) = (1, 1)$, where (see (4.1), (4.2))

$$\varphi_1^A(x, y) = (x-y)^6 = x^6 g_1^A(u), \ g_1^A(u) = (1-u)^6, \ B_1^A = 6,$$
$$\varphi_2^A(x, y) = -(x-y)^2 x^5 = x^7 g_2^A(u), \ g_2^A(u) = -(1-u)^2, \ B_2^A = 7,$$
$$\varphi_3^A(x, y) = x^8 = x^8 g_3^A(u), \ g_3^A(u) = 1, \ B_1^A = 8.$$

The characteristic polynomial $g_1^A(u) = (1-u)^6$ of the quasi-homogeneous form $\varphi_1^A(x, y)$ has a single real root $u_0 = 1$. Note that $A = (A_1, A_2) = (1, 1) \in A_p$, and when checking the condition (C̃4)$_{A,u_0}$, in general we must consider polynomials $x(t)$, $y(t)$ of the form

$$x(t) = c_0 t^v + c_1 t^{v+1} + o(t^{v+1}), \ y(t) = d_0 t^v + d_1 t^{v+1} + o(t^{v+1}), \text{ where } v \in \mathbb{N}. \tag{5.5}$$

Here we have case 3) (when both numbers $e_1 = -A_2 = -1$, $e_2 = A_1 = 1$ are odd), and hence the pair of numbers $(c_0, d_0)$ can be chosen from the set $C_3 = \{(1, 1), (-1, -1)\}$. In case $\nu = 1$ at $c_0 = 1$, $d_0 = 1$, the polynomials (5.5) are of the form (for simplicity we denote by $c_1 = c$, $d_1 = d$)

$$x(t) = t + ct^2 + o(t^2), \quad y(t) = t + dt^2 + o(t^2). \tag{5.6}$$

Then

$$\varphi_1^A(x(t), y(t)) = (x(t) - y(t))^6 = (c-d)^6 t^{12} + o(t^{12}),$$

$$\varphi_2^A(x(t), y(t)) = -(x(t) - y(t))^2 [x(t)]^5 = -(c-d)^2 t^9 + o(t^9),$$

$$\varphi_3^A(x(t), y(t)) = [x(t)]^8 = t^8 + o(t^8).$$

Thus, $p(x(t), y(t)) = (c-d)^6 t^{12} + o(t^{12}) - (c-d)^2 t^9 + o(t^9) + t^8 + o(t^8)$,

and hence for any polynomial of the form (5.6) it is $p(x(t), y(t)) = t^8 + o(t^8)$. Quite similarly, in the case of choosing $c_0 = -1$, $d_0 = -1$, i.e., for polynomials of the form $x(t) = -t + ct^2 + o(t^2)$, $y(t) = -t + dt^2 + o(t^2)$, it is also true that $p(x(t), y(t)) = t^8 + o(t^8)$.

Consider now polynomials of the form (5.5) with $\nu = 2$, $c_0 = 1$, $d_0 = 1$, (again denoted by $c_1 = c$, $d_1 = d$ for simplicity):

$$x(t) = t^2 + ct^3 + o(t^3), \quad y(t) = t^2 + dt^3 + o(t^3). \tag{5.7}$$

Then

$$\varphi_1^A(x(t), y(t)) = (x(t) - y(t))^6 = (c-d)^6 t^{18} + o(t^{18}),$$

$$\varphi_2^A(x(t), y(t)) = -(x(t) - y(t))^2 x^5(t) = -(c-d)^2 t^{16} + o(t^{16}),$$

$$\varphi_3^A(x(t), y(t)) = x^8(t) = t^{16} + o(t^{16}),$$

$$p(x(t), y(t)) = (c-d)^6 t^{18} + o(t^{18}) - (c-d)^2 t^{16} + o(t^{16}) + t^{16} + o(t^{16}),$$

and hence, for example, at $c = 2$, $d = 0$ for polynomials $x(t) = t^2 + 2t^3$, $y(t) = t^2$ is $p(x(t), y(t)) = 2^6 t^{18} - 2^4 t^{16} + t^{16} = -15t^{16} + o(t^{16})$. Thus, the consideration of polynomials of the form (5.7) has shown that $0_{(2)}$ is not a point of local minimum of the polynomial $p(x, y)$, i.e., the consideration of only polynomials of the form (5.6) corresponding to the case $\nu = 1$ was insufficient. Note that when using Algorithm 1 we will not get an answer to the question whether $0_{(2)}$ is the point of local minimum of the polynomial from this example.

**Remark 9.** All the statements used in this section remain valid for the power series $p(x, y)$ [5], and hence can be applied to $p(x, y)$ in this case as well.

## 6. THE CASE WHEN $p(x, y)$ IS THE SUM OF TWO A-QUASI-HOMOGENEOUS FORMS, WHERE $A \in \mathbb{N}_0^2$

Let us now consider one special case when for a polynomial $p(x, y)$, satisfying the conditions of statement 9, for some $A \in \mathbb{N}_0^2$ the decomposition (4.1), (4.2) consists of two A-quasi-homogeneous forms, i.e., $r_A = 2$. Let us show that in this case, by simple computational procedures, we can uniquely answer the question whether $0_{(2)}$ is a point of local minimum of $p(x, y)$.

We will need the following statements.

**Assertion 10.** Let for some $A \in \mathbb{N}_0^2$ the decomposition (4.1), (4.2) of the polynomial $p(x, y)$ consists of two $A$-quasi-homogeneous forms, i.e., has the form $p(x, y) = \varphi_1^A(x, y) + \varphi_2^A(x, y)$. Then there does not exist a vector $\bar{A} \in \mathbb{N}_0^2$ such that $\bar{A} \neq A$, and the main quasi-homogeneous form $\varphi_1^{\bar{A}}(x, y)$ of the polynomial $p(x, y)$ contains more than two single terms.

**Proof.** Suppose that for some $\bar{A} = (\bar{A}_1, \bar{A}_2) \in \mathbb{N}_0^2$ the main $\bar{A}$-quasi-homogeneous form $\varphi_1^{\bar{A}}(x, y)$ of the polynomial $p(x, y)$ contains at least three single terms. At least two of them will belong to one of the two forms $\varphi_1^A(x, y)$ or $\varphi_2^A(x, y)$ (the carrier of each of which lies on the same line). Since a line in the plane is uniquely defined by any two points on this line, the belonging of two single terms to any of these forms means that $\bar{A} = A$ (since $A, \bar{A} \in \mathbb{N}_0^2$), i.e., we have come to a contradiction with $\bar{A} \neq A$.

**Assertion 11.** Let $A = (A_1, A_2) \in \mathbb{N}_0^2$, $u_0 \neq 0$, $l$ be an even natural number, $p(x, y)$, $\bar{p}(x, y)$ be polynomials, $p(x, y) = (y^{A_1} - u_0 x^{A_2})^l \bar{p}(x, y)$, $\bar{p}(0,0) = 0$. Then $0_{(2)}$ is the point of local minimum of the polynomial $p(x, y)$ if and only if it is the point of local minimum of the polynomial $\bar{p}(x, y)$.

**Proof.** 1) Let $0_{(2)}$ not be a point of local minimum of the polynomial $p(x, y)$. Then there exists a sequence of points $(x(n), y(n)) \in \mathbb{R}^2$, such that
$$p(x(n), y(n)) = \left[(y(n))^{A_1} - u_0 (x(n))^{A_2}\right]^l \bar{p}(x(n), y(n)) < 0, \ n = 1, 2, \dots, \ (x(n), y(n)) \to 0_{(2)} \text{ at}$$
$n \to \infty$. From these conditions, using the parity of $l$, we obtain:
$$[y(n)]^{A_1} - u_0 [x(n)]^{A_2} \neq 0, \ \bar{p}((x(n), y(n))) < 0, \ n = 1, 2, \dots,$$
i.e., $0_{(2)}$ is not a point of local minimum of the polynomial $\bar{p}(x, y)$.

2) In the reverse direction. Let $0_{(2)}$ not be a point of local minimum of the polynomial $\bar{p}(x, y)$. Then there is a sequence of points $(x(n), y(n)) \in \mathbb{R}^2$, such that $\bar{p}((x(n), y(n))) < 0$, $n = 1, 2, \dots$, $(x(n), y(n)) \to 0_{(2)}$ at $n \to \infty$. In this case, due to the continuity of $\bar{p}(x, y)$, we can assume that $x(n) \neq 0$, $y(n) \neq 0$, $n = 1, 2, \dots$. Note that if $y_0^{A_1} - u_0 x_0^{A_2} = 0$ is true for some $x_0 \neq 0$, $y_0 \neq 0$, $u_0 \neq 0$, then $\forall v \in (0, 1)$ $(v y_0)^{A_1} - u_0 x_0^{A_2} \neq 0$. Indeed, if $(v y_0)^{A_1} - u_0 x_0^{A_2} = 0$, then $(v y_0)^{A_1} = y_0^{A_1}$, whence $v^{A_1} = 1$ which contradicts the condition $v \in (0, 1)$. But then, using the continuity of $\bar{p}(x, y)$, for any number $n = 1, 2, \dots$, we can find a number $v_n \in (1 - 1/n, 1]$, such that $\bar{p}(x(n), v_n y(n)) < 0$, $(v_n y(n))^{A_1} - u_0 (x(n))^{A_2} \neq 0$ (if $(y(n))^{A_1} - u_0 (x(n))^{A_2} \neq 0$, we assume $v_n = 1$). Thus, we obtain:
$$p(x(n), v_n y(n)) = \left[(v_n y(n))^{A_1} - u_0 (x(n))^{A_2}\right]^l \bar{p}(x(n), v_n y(n)) < 0, \ n = 1, 2, \dots,$$
and $(x(n), y(n)) \to 0_{(2)}$ at $n \to \infty$, i.e., $0_{(2)}$ is not a point of local minimum of the polynomial $p(x, y)$.

Let us give some information concerning an arbitrary $A$-quasi-homogeneous form $\varphi_1^A(x, y)$, where $A = (A_1, A_2) \in \mathbb{N}_0^2$, of the form (2.1) - (2.3). Let, as before, $e = (-A_2, A_1)$, $u = x^{e_1} y^{e_2} = x^{-A_2} y^{A_1}$,

$$\varphi_1^A(x,y) = \sum_{i=1}^s a_i x^{\alpha_i} y^{\beta_i} = x^{\alpha_1} y^{\beta_1} g_1^A(x^{-A_2} y^{A_1}) = x^{\alpha_1} y^{\beta_1} g_1^A(x^{e_1} y^{e_2}) = x^{\alpha_1} y^{\beta_1} g_1^A(u),$$

$$r_1 = \deg g_1^A(u) = v_s = (\alpha_1 - \alpha_s)/A_2, \quad \alpha_1 = \alpha_s + r_1 A_2, \quad \alpha_1 - r_1 A_2 = \alpha_s$$

( $\deg g(u)$ is the degree of an arbitrary polynomial $g(u)$ ). Let, further, $u_0$ be the root of a multiple $k \in \mathbb{N}$ of a polynomial $g_1^A(u)$, i.e., $g_1^A(u) = (u - u_0)^k \bar{g}_1^A(u)$, where $\bar{g}_1^A(u)$ is a polynomial, $\bar{g}_1^A(u_0) \neq 0$, $\bar{r}_1 = \deg \bar{g}_1^A(u) = r_1 - k$. Then

$$\varphi_1^A(x,y) = x^{\alpha_1} y^{\beta_1} g_1^A(x^{-A_2} y^{A_1}) = x^{\alpha_1} y^{\beta_1} (x^{-A_2} y^{A_1} - u_0)^k \bar{g}_1^A(x^{-A_2} y^{A_1}) =$$
$$= (y^{A_1} - u_0 x^{A_2})^k x^{\alpha_1 - r_1 A_2} y^{\beta_1} \left[ x^{\bar{r}_1 A_2} \bar{g}_1^A(x^{-A_2} y^{A_1}) \right] = (y^{A_1} - u_0 x^{A_2})^k \bar{\varphi}_1^A(x,y), \tag{6.1}$$

where the polynomial

$$\bar{\varphi}_1^A(x,y) = x^{\alpha_1 - r_1 A_2} y^{\beta_1} \left[ x^{\bar{r}_1 A_2} \bar{g}_1^A(x^{-A_2} y^{A_1}) \right] = x^{\alpha_1 - k A_2} y^{\beta_1} \bar{g}_1^A(x^{-A_2} y^{A_1}) \tag{6.2}$$

has terms of the form $\bar{a}_i x^{\bar{\alpha}_i} y^{\bar{\beta}_i} = \bar{a}_i x^{\alpha_1 - k A_2} y^{\beta_1} (x^{-A_2} y^{A_1})^{\bar{v}_i} = \bar{a}_i x^{\alpha_1 - (k+\bar{v}_i)A_2} y^{\beta_1 + A_1 \bar{v}_i}$, $\bar{v}_i \in \mathbb{N} \cup \{0\}, \bar{v}_i \leq \bar{r}_1$,

and thus

$$\bar{\alpha}_i = \alpha_1 - (k + \bar{v}_i) A_2 \geq \alpha_1 - (k + \bar{r}_1) A_2 = \alpha_1 - r_1 A_2 = \alpha_s \geq 0, \quad \bar{\beta}_i = \beta_1 + A_1 \bar{v}_i \geq \beta_1,$$
$$B_1^A = A_1 \alpha_1 + A_2 \beta_1 = A_1(\alpha_s + r_1 A_2) + A_2 \beta_1 = r_1 A_1 A_2 + \alpha_s A_1 + A_2 \beta_1,$$
$$A_1 \bar{\alpha}_i + A_2 \bar{\beta}_i = A_1 \left[\alpha_1 - (k+\bar{v}_i)A_2\right] + A_2(\beta_1 + A_1 \bar{v}_i) = A_1\left[\alpha_1 - k A_2\right] + A_2 \beta_1 =$$
$$= A_1 \alpha_1 + A_2 \beta_1 - k A_1 A_2 = B_1^A - k A_1 A_2 = (r_1 - k) A_1 A_2 + \alpha_s A_1 + A_2 \beta_1 \geq 0.$$

Thus, $\bar{\varphi}_1^A(x,y)$ is also an $A$-quasi-homogeneous polynomial form, for the terms $\bar{a}_i x^{\bar{\alpha}_i} y^{\bar{\beta}_i}$ of which at $\bar{B}_1^A = B_1^A - k A_1 A_2$ it follows that

$$A_1 \bar{\alpha}_i + A_2 \bar{\beta}_i = \bar{B}_1^A = (r_1 - k) A_1 A_2 + \alpha_s A_1 + A_2 \beta_1 \geq 0,$$

and $\bar{B}_1^A = 0 \Rightarrow [r_1 = k, \alpha_s = 0, \beta_1 = 0]$. If $\bar{B}_1^A = 0$, then the equality $r_1 = k$ implies that the polynomial $\bar{g}_1^A(u)$ has degree $r_1 - k = 0$, i.e., is a constant $G \neq 0$ (since $\bar{g}_1^A(u_0) \neq 0$), and then

$$\bar{\varphi}_1^A(x,y) = x^{\alpha_1 - k A_2} y^{\beta_1} \bar{g}_1^A(x^{-A_2} y^{A_1}) = G x^{\alpha_1 - k A_2} y^{\beta_1} = G x^{\alpha_s + r_1 A_2 - k A_2} y^{\beta_1} = G x^0 y^0 = G,$$

i.e., the $A$-quasi-homogeneous form $\bar{\varphi}_1^A(x,y)$ is the same constant.

Let $p(x,y)$ be a polynomial for which the conditions of statement 9 are satisfied: $p(x,y) \not\equiv 0$, $p(0,0) = 0$, $p'(0,0) = 0_{(2)}$ (i.e., $0_{(2)}$ is a stationary point), all main quasi-homogeneous forms of the polynomial $p(x,y)$ from groups 1 and 2 are non-negative and nondegenerate in the weak sense, and thus $A_p \neq \emptyset$, $A = (A_1, A_2) \in A_p$. Let, further,

$$p(x,y) = \varphi_1^A(x,y) + \varphi_2^A(x,y) = x^{\alpha_1} y^{\beta_1} g_1^A(u) + x^{\chi_1} y^{\eta_1} g_2^A(u),$$

where $\varphi_1^A(x,y)$ satisfies conditions (4.1) – (4.3) and $\varphi_2^A(x,y)$ satisfies conditions (4.3) – (4.5). As it follows from Assertion 10, in this case $A_p = \{A\}$.

Let $u_0 \in U_p(A)$, and the multiplicity of the root $u_0$ in the polynomial $g_1^A(u)$ is an even natural number $k \in 2\mathbb{N}$ (see Assertion 1), and the multiplicity of the root $u_0$ in the polynomial $g_2^A(u)$ is $l \in \mathbb{N} \cup \{0\}$. Then $g_1^A(u) = (u - u_0)^k \bar{g}_1^A(u)$, $g_2^A(u) = (u - u_0)^l \bar{g}_2^A(u)$, where $\bar{g}_1^A(u_0) > 0$ (by virtue of non-negativity of $g_1^A(u)$; see Assertion 1), $\bar{g}_2^A(u_0) \neq 0$. Let $r_1 = \deg g_1^A(u)$, $\bar{r}_1 = \deg \bar{g}_1^A(u), r_2 = \deg g_2^A(u), \bar{r}_2 = \deg \bar{g}_2^A(u)$. The following 4 cases are possible:

1) $l = 0$, i.e. $g_2^A(u_0) \neq 0$.
2) $0 < l < k$ and $l$ is even.
3) $k \leq l$.
4) $0 < l < k$ and $l$ is odd.

Using the Assertion 11, as well as the above reasoning concerning an arbitrary $A$-quasi-homogeneous form $\varphi_1^A(x, y)$ (see the representation (6.1), (6.2), where $\overline{\varphi}_1^A(x, y)$ is again an $A$-quasi-homogeneous form), it is easy to pass from the polynomial $p(x, y) = \varphi_1^A(x, y) + \varphi_2^A(x, y)$ to the polynomial $\tilde{p}(x, y) = \tilde{\varphi}_1^A(x, y) + \tilde{\varphi}_2^A(x, y)$, which is also the sum of two $A$-quasi-homogeneous forms, such that $U_{\tilde{p}}(A) \subseteq U_p(A)$, the polynomial $\tilde{\varphi}_1^A(x, y)$ is non-negative and for the new polynomial $\tilde{p}(x, y)$ for any $u_0 \in U_{\tilde{p}}(A)$ only cases of the form 1), 4) will be satisfied, and in the case of the form 4) will be $l = 1$. Indeed, using statement 11, for each fixed $u_0 \in U_{\tilde{p}}(A)$, from case 2) we easily pass from case 2) to the case of the form 1); from case 3) to the case where $u_0 \notin U_{\tilde{p}}(A)$; from case 4) at $l > 1$ to the case of the form 4) at $l = 1$. Moreover, by virtue of Assertion 11, under such a transition $0_{(2)}$ is a point of local minimum of the polynomial $p(x, y)$ if and only if it is a point of local minimum of the polynomial $\tilde{p}(x, y)$. If at the same time for a polynomial $\tilde{p}(x, y)$ for any $u_0 \in U_{\tilde{p}}(A)$ the case of the form 1) is satisfied, then we are in the applicability region of Algorithm 1, using which we obtain an unambiguous answer whether $0_{(2)}$ is a point of local minimum of the polynomial $\tilde{p}(x, y)$ (and thus of $p(x, y)$).

Let us now show that if the condition of the form 4) is satisfied for the polynomial $\tilde{p}(x, y)$, where $l = 1$, at least at one $u_0 \in U_{\tilde{p}}(A)$, $0_{(2)}$ is not a point of local minimum of the polynomial $\tilde{p}(x, y)$ (and thus of $p(x, y)$). For simplicity of notation, we assume that $\tilde{p}(x, y) = p(x, y)$.

Thus, we consider the case when (see (2.1) - (2.3), (4.3) - (4.5))
$$p(x, y) = \varphi_1^A(x, y) + \varphi_2^A(x, y) = x^{\alpha_1} y^{\beta_1} g_1^A(u) + x^{\chi_1} y^{\eta_1} g_2^A(u) =$$
$$= x^{\alpha_1} y^{\beta_1} (u - u_0)^k \overline{g}_1^A(u) + x^{\chi_1} y^{\eta_1} (u - u_0) \overline{g}_2^A(u), k \in 2\mathbb{N},$$

where, as before, $e = (-A_2, A_1)$, $u = x^{e_1} y^{e_2} = x^{-A_2} y^{A_1}$, and while using (6.1), (6.2), we have:
$$p(x, y) = (y^{A_1} - u_0 x^{A_2})^k \overline{\varphi}_1^A(x, y) + (y^{A_1} - u_0 x^{A_2}) \overline{\varphi}_2^A(x, y),$$
$$\overline{\varphi}_1^A(x, y) = x^{\alpha_1 - r_1 A_2} y^{\beta_1} \left[ x^{r_1 A_2} \overline{g}_1^A(x^{-A_2} y^{A_1}) \right] = x^{\alpha_1 - kA_2} y^{\beta_1} \overline{g}_1^A(x^{-A_2} y^{A_1}),$$
$$\overline{\varphi}_2^A(x, y) = x^{\chi_1 - r_2 A_2} y^{\eta_1} \left[ x^{(r_2-1)A_2} \overline{g}_2(x^{-A_2} y^{A_1}) \right] = x^{\chi_1 - A_2} y^{\eta_1} \overline{g}_2(x^{-A_2} y^{A_1}).$$

Let $u_0 = c_0^{e_1} d_0^{e_2} = c_0^{-A_2} d_0^{A_1}$ be fulfilled for some $c_0 \neq 0$, $d_0 \neq 0$. Consider the polynomials:
$$x(t) = c_0 t^{A_1}, y(t) = d_0 t^{A_2} \pm d_0 t^{A_2 + \kappa} = d_0 t^{A_2} (1 \pm t^\kappa), \kappa \in \mathbb{N}.$$

Then (specific values of some quantities: $h_1, \ldots, h_5, \sigma_1, \ldots, \sigma_5$ do not matter)
$$[x(t)]^{-A_2} = [c_0 t^{A_1}]^{-A_2} = c_0^{-A_2} t^{-A_1 A_2}, [y(t)]^{A_1} = [d_0 t^{A_2}(1 \pm t^\kappa)]^{A_1} = d_0^{A_1} t^{A_1 A_2} + o(t^{A_1 A_2}),$$
$$[x(t)]^{-A_2} [y(t)]^{A_1} = c_0^{-A_2} d_0^{A_1} + O(t) = u_0 + O(t),$$
$$\overline{g}_1^A \left( [x(t)]^{-A_2} [y(t)]^{A_1} \right) = \overline{g}_1^A (u_0 + O(t)) = \overline{g}_1^A(u_0) + O(t),$$
$$\overline{g}_2^A \left( [x(t)]^{-A_2} [y(t)]^{A_1} \right) = \overline{g}_2^A (u_0 + O(t)) = \overline{g}_2^A(u_0) + O(t),$$

$$\bar{\varphi}_1^A(x(t), y(t)) = [x(t)]^{\alpha_1 - kA_2}[y(t)]^{\beta_1} \bar{g}_1^A\left([x(t)]^{-A_2}[y(t)]^{A_1}\right) = h_1 t^{\sigma_1} + o(t^{\sigma_1}),$$

$h_1 > 0$, $\sigma_1 \in \mathbb{N} \cup \{0\}$,

$$\bar{\varphi}_2^A(x(t), y(t)) = [x(t)]^{\chi_1 - A_2}[y(t)]^{\eta_1} \bar{g}_2^A\left([x(t)]^{-A_2}[y(t)]^{A_1}\right) = h_2 t^{\sigma_2} + o(t^{\sigma_2}),$$

$h_2 \neq 0$, $\sigma_2 \in \mathbb{N} \cup \{0\}$,

$$[y(t)]^{A_1} - u_0[x(t)]^{A_2} = d_0^{A_1} t^{A_1 A_2}(1 \pm t^\kappa)^{A_1} - u_0 c_0^{A_2} t^{A_1 A_2} =$$
$$= d_0^{A_1} t^{A_1 A_2}(1 \pm A_1 t^\kappa + o(t^\kappa)) - u_0 c_0^{A_2} t^{A_1 A_2} =$$
$$= \left[d_0^{A_1} - u_0 c_0^{A_2}\right] t^{A_1 A_2} \pm d_0^{A_1} A_1 t^{A_1 A_2 + \kappa} + o(t^{A_1 A_2 + \kappa}) = \pm d_0^{A_1} A_1 t^{A_1 A_2 + \kappa} + o(t^{A_1 A_2 + \kappa}).$$

Consequently,

$$\left([y(t)]^{A_1} - u_0[x(t)]^{A_2}\right)^k = h_3 t^{\sigma_3 + k\kappa} + o(t^{\sigma_3 + k\kappa}), \quad h_3 > 0, \quad \sigma_3 \in \mathbb{N}.$$

Thus, we obtain

$$p(x(t), y(t)) = \left([y(t)]^{A_1} - u_0[x(t)]^{A_2}\right)^k \bar{\varphi}_1^A(x(t), y(t)) +$$
$$+ \left([y(t)]^{A_1} - u_0[x(t)]^{A_2}\right) \bar{\varphi}_2^A(x(t), y(t)) = \qquad (6.3)$$
$$= h_4 t^{\sigma_4 + k\kappa} + h_5 t^{\sigma_5 + \kappa} + o(t^{\sigma_4 + k\kappa}) + o(t^{\sigma_5 + \kappa}),$$

where at an appropriate choice of sign in the polynomial $y(t) = d_0 t^{A_2}(1 \pm t^\kappa)$ (we take into account the sign of $h_2 \neq 0$), $h_4 > 0$, $h_5 < 0$, $\sigma_4, \sigma_5 \in \mathbb{N}$ is satisfied. Let us choose a number $\kappa \in \mathbb{N}$ so large that $\sigma_4 + k\kappa > \sigma_5 + \kappa$. This can be done since $k \geq 2$ (i.e., it suffices to take any $\kappa > \sigma_5 - \sigma_4$). Then from (6.3) we obtain

$$p(x(t), y(t)) = h_5 t^{\sigma_5 + \kappa} + o(t^{\sigma_5 + \kappa}), \quad h_5 < 0, \quad \sigma_5 + \kappa \in \mathbb{N},$$

i.e., condition (C4) is satisfied, and hence $0_{(2)}$ is not a point of local minimum of the polynomial $p(x, y)$.